\documentclass[12pt]{article}
\usepackage{amsmath,amssymb}
\usepackage{graphicx}
\newcommand{\eqdef}{\stackrel{\text{def}}{=}}

\newcommand{\n}{\nonumber\\}
\newcommand{\bm}{\boldsymbol}
\newcommand{\cF}{c_{\text{\tiny$\mathcal{F}$}}}
\newcommand{\raprod}[2]{\prod_{#1}^{\stackrel{#2}{\longrightarrow}}}

\newcommand{\ignore}[1]{}
\numberwithin{equation}{section}
\newcommand{\Romannumeral}[1]{\uppercase\expandafter{\romannumeral#1}}
\newcommand{\I}{\text{\Romannumeral{1}}}
\newcommand{\II}{\text{\Romannumeral{2}}}

\newtheorem{thm}{\bf Theorem}
\newcounter{mybangouI}
\newcounter{mythmbangou}

\newtheorem{mythm}{\bf Theorem}[mybangouI]
\newcounter{mybangouII}
\newcounter{myrembangou}
\renewcommand{\themyrembangou}{\arabic{mybangouII}.\arabic{myrembangou}}
\newcommand{\remark}{\refstepcounter{myrembangou}
\noindent{\bf Remark\,\,\themyrembangou}\ \ }

\allowdisplaybreaks[4]

\begin{document}

\baselineskip=20pt

%%%%%%%%%%%%%%%%%%%%%%%%%%%%%%%%%%%%%%%%%%%%%%%%%%%%%%%%%%%%
%                                                          %
%  Title page                                              %
%                                                          %
%%%%%%%%%%%%%%%%%%%%%%%%%%%%%%%%%%%%%%%%%%%%%%%%%%%%%%%%%%%%
\newcommand{\preprint}{
\vspace*{-20mm}
   \begin{flushright}\normalsize \sf
    DPSU-22-3\\
%     {\tt arXiv:2301.00678[math.CA]}\\
%     December 2022
  \end{flushright}}
\newcommand{\Title}[1]{{\baselineskip=26pt
  \begin{center} \Large \bf #1 \\ \ \\ \end{center}}}
\newcommand{\Author}{\begin{center}
  \large \bf Satoru Odake \end{center}}
\newcommand{\Address}{\begin{center}
%     Department of Physics, Shinshu University,\\
     Faculty of Science, Shinshu University,
     Matsumoto 390-8621, Japan
   \end{center}}
\newcommand{\Accepted}[1]{\begin{center}
  {\large \sf #1}\\ \vspace{1mm}{\small \sf Accepted for Publication}
  \end{center}}

\preprint
\thispagestyle{empty}

%\Title{Another Type of Forward and Backward Shift Relations for Orthogonal
%Polynomials in the Askey Scheme}
\Title{Another Type of Forward and Backward\\ Shift Relations for Orthogonal
Polynomials\\ in the Askey Scheme}

\Author

\Address
\vspace{1cm}

\begin{abstract}
The forward and backward shift relations are basic properties of the
(basic) hypergeometric orthogonal polynomials in the Askey scheme (Jacobi,
Askey-Wilson, $q$-Racah, big $q$-Jacobi etc.) and they are related
to the factorization of the differential or difference operators.
Based on other factorizations, we obtain another type of forward and
backward shift relations.
Essentially, these shift relations shift only the parameters.
\end{abstract}

%%%%%%%%%%%%%%%%%%%%%%%%%%%%%%%%%%%%%%%%%%%%%%%%%%%%%%%%%%%%%%%
%                                                             %
%  1. Introduction                                            %
%                                                             %
%%%%%%%%%%%%%%%%%%%%%%%%%%%%%%%%%%%%%%%%%%%%%%%%%%%%%%%%%%%%%%%
\section{Introduction}
\label{sec:intro}

The (basic) hypergeometric orthogonal polynomials in the Askey scheme satisfy
second order differential or difference equations and the forward and backward
shift relations are their basic properties \cite{ismail,kls}.
The orthogonal polynomials in the Askey scheme provide us with exactly solvable
quantum mechanical models. Conversely, we can use the quantum mechanical
formulation as a tool to investigate orthogonal polynomials \cite{os24}.
For example, the forward and backward shift relations are a consequence of the
%shape invariance \cite{os12,os13,os24}, and the multi-indexed orthogonal
%polynomials (\cite{os25,os26,os27} etc.) are found by using the quantum
%mechanical formulation.
shape invariance, and the multi-indexed orthogonal polynomials are found by
using the quantum mechanical formulation.
The Schr\"odinger equation is a second order differential equation for
ordinary quantum mechanics (oQM) and a second order difference equation for
discrete quantum mechanics (dQM). There are two types of dQM, dQM with pure
imaginary shifts (idQM) and dQM with real shifts (rdQM) \cite{os24}.
The coordinate $x$ for oQM and idQM is continuous and that for rdQM is discrete.

The forward and backward shift relations are related to the factorization of
the Hamiltonian. Recently another factorization of the Hamiltonian was found
in a study of the state-adding Darboux transformations for the finite rdQM
systems \cite{os40}. It gives another forward and backward shift relations for
the orthogonal polynomials appearing in the finite rdQM systems ($q$-Racah
etc.), which were called the forward and backward $x$-shift
relations \cite{os40}.
In this paper, we investigate whether such new factorization and forward and
backward shift relations exist for other orthogonal polynomials.
In addition to the finite rdQM systems ($q$-Racah etc.), we examine the oQM
systems (Jacobi etc.), the idQM systems (Askey-Wilson etc.), the semi-infinite
rdQM systems ($q$-Meixner etc.) and the rdQM systems with the Jackson integral
type measure (big $q$-Jacobi etc.).
We call the last category rdQMJ. The quantum mechanical formulation of the
rdQMJ systems needs two component formalism \cite{os34}.
We consider all the polynomials in chapter 9 and 14 of \cite{kls} and
the dual quantum $q$-Krawtchouk polynomial.

This paper is organized as follows.
The orthogonal polynomials in the Askey scheme and their second order
differential or difference equations are recalled in section \ref{sec:OP}.
The forward and backward shift relations are reviewed in section \ref{sec:fbsr}.
Section \ref{sec:newfbsr} is the main part of this paper and new factorization
and another type of forward and backward shift relations are presented.
Section \ref{sec:summary} is for a summary and comments.
In Appendix\,\ref{app:data} the data for \S\,\ref{sec:newfbsr} are given.

%%%%%%%%%%%%%%%%%%%%%%%%%%%%%%%%%%%%%%%%%%%%%%%%%%%%%%%%%%%%%%%
%                                                             %
%  2. Orthogonal Polynomials in the Askey Scheme              %
%                                                             %
%%%%%%%%%%%%%%%%%%%%%%%%%%%%%%%%%%%%%%%%%%%%%%%%%%%%%%%%%%%%%%%
\section{Orthogonal Polynomials in the Askey Scheme}
\label{sec:OP}

In this section we fix the notation and recall the second order differential
or difference equations for the orthogonal polynomials in the Askey scheme
\cite{ismail,kls}.

In our quantum mechanical formulation \cite{os24}, the orthogonal polynomials
in the Askey scheme are expressed as
\begin{equation}
  \check{P}_n(x;\bm{\lambda})\eqdef
  P_n\bigl(\eta(x;\bm{\lambda});\bm{\lambda}\bigr)
  \text{ : a polynomial of degree $n$ in $\eta(x;\bm{\lambda})$}
  \label{cPn=Pn}
\end{equation}
for $n\in\mathbb{Z}_{\geq 0}$ and $\check{P}_{-1}(x;\bm{\lambda})\eqdef
P_{-1}(\eta(x;\bm{\lambda});\bm{\lambda})\eqdef 0$.
Here $x$ is a coordinate of quantum mechanical system and $\eta(x)$ is a
sinusoidal coordinate \cite{os7}, and
$\bm{\lambda}=(\lambda_1,\lambda_2,\ldots)$ are parameters, whose dependence
is expressed as $f=f(\bm{\lambda})$ and $f(x)=f(x;\bm{\lambda})$.
The parameter $q$ is $0<q<1$ and $q^{\bm{\lambda}}$ stands for
$q^{(\lambda_1,\lambda_2,\ldots)}=(q^{\lambda_1},q^{\lambda_2},\ldots)$, and
we omit writing $q$-dependence.

We consider the following orthogonal polynomials, all the polynomials in
chapter 9 and 14 of \cite{kls} and the dual quantum $q$-Krawtchouk polynomial:
Hermite (He),
Laguerre (L),
Jacobi (J),
Bessel (B),
pseudo Jacobi (pJ),
continuous Hahn (cH),
Meixner-Pollaczek (MP),
Wilson (W),
continuous dual Hahn (cdH),
Askey-Wilson (AW),
continuous dual $q$-Hahn (cd$q$H),
Al-Salam-Chihara (ASC),
continuous big $q$-Hermite (cb$q$He),
continuous $q$-Hermite (c$q$He),
continuous $q$-Jacobi (c$q$J),
continuous $q$-Laguerre (c$q$L),
continuous $q$-Hahn (c$q$H),
$q$-Meixner-Pollaczek ($q$MP),
Hahn (H),
Krawtchouk (K),
Racah (R),
dual Hahn (dH),
dual quantum $q$-Krawtchouk (dq$q$K) (which is not treated in \cite{kls}),
$q$-Hahn ($q$H),
$q$-Krawtchouk ($q$K),
quantum $q$-Krawtchouk (q$q$K),
affine $q$-Krawtchouk (a$q$K),
$q$-Racah ($q$R),
dual $q$-Hahn (d$q$H),
dual $q$-Krawtchouk (d$q$K),
Meixner (M),
Charlier (C),
little $q$-Jacobi (l$q$J),
little $q$-Laguerre/Wall (l$q$L),
$q$-Bessel ($q$B) (=alternative $q$-Charlier),
$q$-Meixner ($q$M),
Al-Salam-Carlitz $\II$ (ASC$\II$),
$q$-Charlier ($q$C),
big $q$-Jacobi (b$q$J),
big $q$-Laguerre (b$q$L),
Al-Salam-Carlitz $\I$ (ASC$\I$),
discrete $q$-Hermite $\I$ (d$q$He$\I$),
discrete $q$-Hermite $\II$ (d$q$He$\II$),
$q$-Laguerre ($q$L)
and Stieltjes-Wigert (SW).
Explicit expressions of various quantities ($\check{P}_n(x)$, $P_n(\eta)$,
$\eta(x)$, $\varphi(x)$, $\mathcal{E}_n$, $\bm{\lambda}$, $\bm{\delta}$,
$\kappa$, $f_n$, $b_n$, $c_1(\eta)$, $c_2(\eta)$, $\cF$, $w(x)$, $V(x)$, $B(x)$,
$D(x)$, $B^{\text{J}}(x)$, $D^{\text{J}}(x)$, $f^{\text{J}}_n$,
$b^{\text{J}}_n$) are given in Appendix A of arXiv:2301.00678v1.
These polynomials appear in quantum mechanical systems as follows:
\begin{align}
  \text{oQM}&:\ \text{He,\,L,\,J,\,B,\,pJ},
  \label{polyoQM}\\
  \text{idQM}&:\ \text{cH,\,MP,\,W,\,cdH,\,AW,\,cd$q$H,\,ASC,\,cb$q$He,
  c$q$He,\,c$q$J,\,c$q$L,\,c$q$H,\,$q$MP},
  \label{polyidQM}\\
  \text{rdQM}&\text{ (finite)}:\ \text{H,\,K,\,R,\,dH,\,dq$q$K,\,$q$H,\,$q$K,
  q$q$K,\,a$q$K,\,$q$R,\,d$q$H,\,d$q$K},
  \label{polyrdQMf}\\
  \text{rdQM}&\text{ (semi-infinite)}:\ \text{M,\,C,\,l$q$J,\,l$q$L,\,$q$B,
  $q$M,\,ASC$\II$,\,$q$C},
  \label{polyrdQMsi}\\
  \text{rdQMJ}&:\ \text{b$q$J,\,b$q$L,\,ASC$\I$,\,d$q$He$\I$,\,d$q$He$\II$,
  $q$L,\,SW}.
  \label{polyrdQMJ}
\end{align}
We comment that the oQM systems described by the Bessel and pseudo Jacobi
polynomials are the Morse potential and the hyperbolic symmetric top $\II$,
respectively.
We also comment on an infinite sum orthogonality relations for
the Stieltjes-Wigert polynomial (parameter: $c>0$),
\begin{equation}
  \text{SW}:
  \ \sum_{x=-\infty}^{\infty}c^xq^{\frac12x(x+1)}P_n(cq^x)P_m(cq^x)
  =\delta_{nm}\,q^{-n}(q\,;q)_n(q,-cq,-c^{-1}\,;q)_{\infty}
  \ \ (n,m\in\mathbb{Z}_{\geq 0}),
  \label{orthoSW}
\end{equation}
which are obtained from those for $q$-Laguerre polynomial by taking an
appropriate limit. This \eqref{orthoSW} is not given in \cite{kls}.

The Schr\"odinger equations of oQM and dQM systems are second order
differential and difference equations, respectively.
By the similarity transformation in terms of the ground state wavefunction,
the similarity transformed Hamiltonian $\widetilde{\mathcal{H}}(\bm{\lambda})$
is a second order differential or difference operator acting on the
eigenpolynomials $\check{P}_n(x;\bm{\lambda})$ \cite{os24},
\begin{alignat}{2}
  \text{oQM}&:\ \ &
  \widetilde{\mathcal{H}}(\bm{\lambda})&\eqdef 
  -\frac{d^2}{dx^2}-2\frac{dw(x;\bm{\lambda})}{dx}\frac{d}{dx}
  \ \Bigl(=-4c_2(\eta)\frac{d^2}{d\eta^2}
  -4c_1(\eta;\bm{\lambda})\frac{d}{d\eta}\Bigr),
  \label{tHoQM}\\
  \text{idQM}&:\ \ &
  \widetilde{\mathcal{H}}(\bm{\lambda})&\eqdef
  V(x;\bm{\lambda})(e^{\gamma p}-1)+V^*(x;\bm{\lambda})(e^{-\gamma p}-1),
  \label{tHidQM}\\
  \text{rdQM}&:\ \ &
  \widetilde{\mathcal{H}}(\bm{\lambda})&\eqdef
  B(x;\bm{\lambda})(1-e^{\partial})+D(x;\bm{\lambda})(1-e^{-\partial}).
  \label{tHrdQM}
\end{alignat}
For oQM, the coordinate $x$ is a continuous variable and the Hamiltonian
$\mathcal{H}(\bm{\lambda})$ is
\begin{equation}
  \mathcal{H}(\bm{\lambda})=-\frac{d^2}{dx^2}+U(x;\bm{\lambda}),\quad
  U(x;\bm{\lambda})\eqdef\Bigl(\frac{dw(x;\bm{\lambda})}{dx}\Bigr)^2
  +\frac{d^2w(x;\bm{\lambda})}{dx^2}.
  \label{Udef}
\end{equation}
While the orthogonality relations of $\check{P}_n(x)$ for B and pJ cases hold
only for a finite number of $n$, we consider all
$n\in\mathbb{Z}_{\geq 0}$, because we consider only differential equations
(or relations) in this paper.
For idQM, the coordinate $x$ is a continuous variable and the momentum $p$ is
$p=-i\frac{d}{dx}$, and $\gamma$ is a real constant ($\gamma=1$ for non
$q$-polynomial, $\gamma=\log q$ for $q$-polynomial). The operator $e^{\alpha p}$
($\alpha$: constant) is a shift operator, $e^{\alpha p}f(x)=f(x-i\alpha)$.
The $*$-operation on an analytic function $f(x)=\sum_na_nx^n$
($a_n\in\mathbb{C}$) is defined by $f^*(x)=\sum_na_n^*x^n$, in which
$a_n^*$ is the complex conjugation of $a_n$.
For rdQM, the Schr\"odinger equation is a matrix eigenvalue problem.
The similarity transformed Hamiltonian
$\widetilde{\mathcal{H}}=(\widetilde{\mathcal{H}}_{x,y})$ is a matrix labeled
by the coordinate $x$, which takes discrete values in $\{0,1,\ldots,N\}$ or
$\mathbb{Z}_{\geq 0}$. In this paper, however, we treat $x$ as a continuous
variable $x\in\mathbb{R}$, because we only deal with difference equations
(or relations). The operators $e^{\pm\partial}$ are shift operators
$e^{\pm\partial}=e^{\pm\frac{d}{dx}}$, $e^{\pm\partial}f(x)=f(x\pm 1)$.
We consider $\check{P}_n(x)$ with all $n\in\mathbb{Z}_{\geq 0}$ even for
finite systems.
We remark that the polynomials $\check{P}_n(x)$ in finite rdQM \eqref{polyrdQMf},
whose orthogonality holds for $n=0,1,\ldots,N$, are
ill-defined for $n>N$ due to the normalization condition
$\check{P}_n(0)=P_n(0)=1$.
So we should replace $\check{P}_n(x)$ ($n>N$) in finite rdQM
with the monic version
$\check{P}^{\text{monic}}_n(x;\bm{\lambda})\eqdef
c_n(\bm{\lambda})^{-1}\check{P}_n(x;\bm{\lambda})$
($c_n(\bm{\lambda})$: the coefficient of the highest degree term)
in Theorem\,\ref{thm:sabuneq}, \ref{thm:fbsr} and \ref{thm:newfbsr_rdQM}
(with the replacements $f_n(\bm{\lambda})\to f^{\text{monic}}_n(\bm{\lambda})
=f_n(\bm{\lambda})c_n(\bm{\lambda})^{-1}c_{n-1}(\bm{\lambda}+\bm{\delta})$,
etc.).

For rdQMJ (the rdQM system with Jackson integral type measure such as the
big $q$-Jacobi polynomial), its quantum mechanical formulation needs
two component formalism with two sinusoidal coordinates
$\eta^{(\pm)}(x;\bm{\lambda})$ \cite{os34}.
Since only difference equations (or relations) are considered in this paper,
we use $\eta$ only (we do not use $x$) and treat $\eta$ as a continuous
variable $\eta\in\mathbb{R}$.
The similarity transformed Hamiltonian
$\widetilde{\mathcal{H}}^{\text{J}}(\bm{\lambda})$ is a second order difference
operator acting on the eigenpolynomials $P_n(\eta;\bm{\lambda})$ \cite{os34},
\begin{equation}
  \text{rdQMJ}:\ \ \widetilde{\mathcal{H}}^{\text{J}}(\bm{\lambda})
  \eqdef B^{\text{J}}(\eta;\bm{\lambda})(1-q^{\eta\frac{d}{d\eta}})
  +D^{\text{J}}(\eta;\bm{\lambda})(1-q^{-\eta\frac{d}{d\eta}}),
  \label{tHrdQMJ}
\end{equation}
where the operators $q^{\pm\eta\frac{d}{d\eta}}$ are $q$-shift operators,
$q^{\pm\eta\frac{d}{d\eta}}f(\eta)=f(q^{\pm1}\eta)$.

The orthogonal polynomials in the Askey scheme studied in this paper have the
following property.
\begin{thm}\cite{ismail,kls}
\label{thm:sabuneq}
The polynomials in \eqref{polyoQM}--\eqref{polyrdQMJ} satisfy the second order
differential or difference equations for $n\in\mathbb{Z}_{\geq 0}$,
\begin{align}
  \text{\rm oQM,\,idQM,\,rdQM}&:
  \ \widetilde{\mathcal{H}}(\bm{\lambda})\check{P}_n(x;\bm{\lambda})
  =\mathcal{E}_n(\bm{\lambda})\check{P}_n(x;\bm{\lambda}),
  \label{tHcPn=}\\
  \text{\rm rdQMJ}&:
  \widetilde{\mathcal{H}}^{\text{\rm J}}(\bm{\lambda})P_n(\eta;\bm{\lambda})
  =\mathcal{E}_n(\bm{\lambda})P_n(\eta;\bm{\lambda}).
  \label{tHJPn=}
\end{align}
\end{thm}
We remark that the constant terms of $\widetilde{\mathcal{H}}$ and
$\widetilde{\mathcal{H}}^{\text{\rm J}}$ are chosen such that $\mathcal{E}_0=0$.
For idQM, the relation \eqref{tHcPn=} is invariant under the $*$-operation.

%%%%%%%%%%%%%%%%%%%%%%%%%%%%%%%%%%%%%%%%%%%%%%%%%%%%%%%%%%%%%%%
%                                                             %
%  3. Forward and Backward Shift Relations                    %
%                                                             %
%%%%%%%%%%%%%%%%%%%%%%%%%%%%%%%%%%%%%%%%%%%%%%%%%%%%%%%%%%%%%%%
\section{Forward and Backward Shift Relations}
\label{sec:fbsr}
\setcounter{mybangouI}{2}
\setcounter{mythm}{0}

The similarity transformed Hamiltonians $\widetilde{\mathcal{H}}(\bm{\lambda})$
\eqref{tHoQM}--\eqref{tHrdQM} are factorized as
\begin{equation}
  \widetilde{\mathcal{H}}(\bm{\lambda})
  =\mathcal{B}(\bm{\lambda})\mathcal{F}(\bm{\lambda}),
  \label{tH=BF}
\end{equation}
where the forward and backward shift operators, $\mathcal{F}(\bm{\lambda})$ and
$\mathcal{B}(\bm{\lambda})$, are defined by \cite{os12,os13},
\begin{alignat}{2}
  \text{oQM}&:\ \ &
  \mathcal{F}(\bm{\lambda})&\eqdef
  \cF\Bigl(\frac{d\eta(x)}{dx}\Bigr)^{-1}\frac{d}{dx}
  \ \Bigl(=\cF\frac{d}{d\eta}\Bigr),
  \label{FoQM}\\
  &&\mathcal{B}(\bm{\lambda})&\eqdef
  -\cF^{-1}\Bigl(\frac{d\eta(x)}{dx}\frac{d}{dx}
  +4c_1\bigl(\eta(x);\bm{\lambda}\bigr)\Bigr)
  \ \Bigl(=-4\cF^{-1}\Bigl(c_2(\eta)\frac{d}{d\eta}
  +c_1(\eta;\bm{\lambda})\Bigr)\Bigr),
  \label{BoQM}\\
  \text{idQM}&:\ \ &
  \mathcal{F}(\bm{\lambda})&\eqdef i\varphi(x;\bm{\lambda})^{-1}
  (e^{\frac{\gamma}{2}p}-e^{-\frac{\gamma}{2}p}),
  \label{FidQM}\\
  &&\mathcal{B}(\bm{\lambda})&\eqdef -i
  \bigl(V(x;\bm{\lambda})e^{\frac{\gamma}{2}p}
  -V^*(x;\bm{\lambda})e^{-\frac{\gamma}{2}p}\bigr)\varphi(x;\bm{\lambda}),
  \label{BidQM}\\
  \text{rdQM}&:\ \ &
  \mathcal{F}(\bm{\lambda})&\eqdef
  B(0;\bm{\lambda})\varphi(x;\bm{\lambda})^{-1}(1-e^{\partial}),
  \label{FrdQM}\\
  &&\mathcal{B}(\bm{\lambda})&\eqdef\frac{1}{B(0;\bm{\lambda})}
  \bigl(B(x;\bm{\lambda})-D(x;\bm{\lambda})e^{-\partial}\bigr)
  \varphi(x;\bm{\lambda}).
  \label{BrdQM}
\end{alignat}
Since $w(x)$, $B(x)$, $D(x)$ and $V(x)$ satisfy
\begin{align}
  &\Bigl(\frac{dw(x;\bm{\lambda})}{dx}\Bigr)^2
  -\frac{d^2w(x;\bm{\lambda})}{dx^2}
  =\Bigl(\frac{dw(x;\bm{\lambda}+\bm{\delta})}{dx}\Bigr)^2
  +\frac{d^2w(x;\bm{\lambda}+\bm{\delta})}{dx^2}
  +\mathcal{E}_1(\bm{\lambda}),\\
  &\frac{V(x-i\frac{\gamma}{2};\bm{\lambda})}{V(x;\bm{\lambda}+\bm{\delta})}
  =\kappa\frac{\varphi(x;\bm{\lambda})}{\varphi(x-i\gamma;\bm{\lambda})},\\
  &V(x+i\tfrac{\gamma}{2};\bm{\lambda})+V^*(x-i\tfrac{\gamma}{2};\bm{\lambda})
  =\kappa\bigl(V(x;\bm{\lambda}+\bm{\delta})
  +V^*(x;\bm{\lambda}+\bm{\delta})\bigr)-\mathcal{E}_1(\bm{\lambda}),\\
  &\frac{B(x+1;\bm{\lambda})}{B(x;\bm{\lambda}+\bm{\delta})}
  =\kappa\frac{\varphi(x;\bm{\lambda})}{\varphi(x+1;\bm{\lambda})},\quad
  \frac{D(x;\bm{\lambda})}{D(x;\bm{\lambda}+\bm{\delta})}
  =\kappa\frac{\varphi(x;\bm{\lambda})}{\varphi(x-1;\bm{\lambda})},\\
  &B(x;\bm{\lambda})+D(x+1;\bm{\lambda})=
  \kappa\bigl(B(x;\bm{\lambda}+\bm{\delta})+D(x;\bm{\lambda}+\bm{\delta})\bigr)
  +\mathcal{E}_1(\bm{\lambda}),
\end{align}
we obtain
\begin{equation}
  \mathcal{F}(\bm{\lambda})\mathcal{B}(\bm{\lambda})
  =\kappa\mathcal{B}(\bm{\lambda}+\bm{\delta})
  \mathcal{F}(\bm{\lambda}+\bm{\delta})
  +\mathcal{E}_1(\bm{\lambda}),
  \label{shapeinv}
\end{equation}
which is the (similarity transformed) shape invariance condition.
The constants $f_n$ and $b_n$ satisfy
\begin{equation}
  \mathcal{E}_n(\bm{\lambda})=f_n(\bm{\lambda})b_{n-1}(\bm{\lambda})
  \ \ (n\in\mathbb{Z}_{\geq 0}),
\end{equation}
and the energy eigenvalues $\mathcal{E}_n$ satisfy
\begin{equation}
  \mathcal{E}_{n+1}(\bm{\lambda})=\kappa\mathcal{E}_n(\bm{\lambda}+\bm{\delta})
  +\mathcal{E}_1(\bm{\lambda})\ \ (n\in\mathbb{Z}_{\geq 0}).
  \label{Enshapeinv}
\end{equation}
Note that we have $f_n(\bm{\lambda})=\mathcal{E}_n(\bm{\lambda})$ and
$b_n(\bm{\lambda})=1$ for rdQM due to our normalization
$\check{P}_n(0)=P_n(0)=1$.
Corresponding to the factorization \eqref{tH=BF}, the shape invariance combined
with the Crum's theorem give the following relations \cite{os12,os13,os24}.
\begin{mythm}\cite{kls}
\label{thm:fbsr}
For the polynomials in \eqref{polyoQM}--\eqref{polyrdQMsi},
the forward and backward shift relations hold:
\begin{align}
  \mathcal{F}(\bm{\lambda})\check{P}_n(x;\bm{\lambda})
  &=f_n(\bm{\lambda})\check{P}_{n-1}(x;\bm{\lambda}+\bm{\delta})
  \ \ (n\in\mathbb{Z}_{\geq 0}),
  \label{FcPn=}\\
  \mathcal{B}(\bm{\lambda})\check{P}_{n-1}(x;\bm{\lambda}+\bm{\delta})
  &=b_{n-1}(\bm{\lambda})\check{P}_n(x;\bm{\lambda})
  \ \ (n\in\mathbb{Z}_{\geq 1}).
  \label{BcPn-1=}
\end{align}
\end{mythm}
For idQM, the relations \eqref{FcPn=}--\eqref{BcPn-1=} are invariant under the
$*$-operation.

\medskip

The similarity transformed Hamiltonians
$\widetilde{\mathcal{H}}^{\text{J}}(\bm{\lambda})$ \eqref{tHrdQMJ} are
factorized as
\begin{equation}
  \widetilde{\mathcal{H}}^{\text{J}}(\bm{\lambda})
  =\mathcal{B}^{\text{J}}(\bm{\lambda})\mathcal{F}^{\text{J}}(\bm{\lambda}).
  \label{tHJ=BJFJ}
\end{equation}
Here the forward and backward shift operators,
$\mathcal{F}^{\text{J}}(\bm{\lambda})$ and
$\mathcal{B}^{\text{J}}(\bm{\lambda})$, are defined by \cite{os34}
\begin{alignat}{2}
  \text{rdQMJ}&:\ \ &
  \mathcal{F}^{\text{J}}(\bm{\lambda})&\eqdef
  A(\bm{\lambda})\eta^{-1}(1-q^{\eta\frac{d}{d\eta}}),
  \label{FJrdQMJ}\\
  &&\mathcal{B}^{\text{J}}(\bm{\lambda})&\eqdef
  \frac{1}{A(\bm{\lambda})}\bigl(
  B^{\text{J}}(\eta;\bm{\lambda})-D^{\text{J}}(\eta;\bm{\lambda})
  q^{-\eta\frac{d}{d\eta}}\bigr)\eta,
  \label{BJrdQMJ}
\end{alignat}
where the constant $A(\bm{\lambda})$ is given by 
\begin{equation}
  A(\bm{\lambda})=\left\{
  \begin{array}{ll}
  -q^{-1}D^{\text{J}}(1;\bm{\lambda})&:\text{b$q$J,\,b$q$L,\,$q$L,\,SW}\\
  -aq^{-1}&:\text{ASC$\I$}\\
  q^{-1}&:\text{d$q$He$\I$}\\
  1&:\text{d$q$He$\II$}
 \end{array}\right..
\end{equation}
We can show that $B^{\text{J}}(\eta)$ and $D^{\text{J}}(\eta)$ satisfy
\begin{align}
  &qB^{\text{J}}(qr^{-1}\eta;\bm{\lambda})
  =\kappa B^{\text{J}}(\eta;\bm{\lambda}+\bm{\delta}),\quad
  q^{-1}D^{\text{J}}(r^{-1}\eta;\bm{\lambda})
  =\kappa D^{\text{J}}(\eta;\bm{\lambda}+\bm{\delta}),\\
  &B^{\text{J}}(r^{-1}\eta;\bm{\lambda})+D^{\text{J}}(qr^{-1}\eta;\bm{\lambda})
  =\kappa\bigl(B^{\text{J}}(\eta;\bm{\lambda}+\bm{\delta})
  +D^{\text{J}}(\eta;\bm{\lambda}+\bm{\delta})\bigr)
  +\mathcal{E}_1(\bm{\lambda}),
\end{align}
where $r$ is given by
\begin{equation}
  r=\left\{
  \begin{array}{ll}
  q&:\text{b$q$J,\,b$q$L,\,d$q$He$\II$,\,$q$L}\\[2pt]
  1&:\text{ASC$\I$,\,d$q$He$\I$}\\[2pt]
  q^2&:\text{SW}
  \end{array}\right..
\end{equation}
Therefore we obtain
\begin{equation}
  \mathcal{F}^{\text{J}}(\bm{\lambda})\mathcal{B}^{\text{J}}(\bm{\lambda})
  \Bigm|_{\eta\to r^{-1}\eta}
  =\kappa\mathcal{B}^{\text{J}}(\bm{\lambda}+\bm{\delta})
  \mathcal{F}^{\text{J}}(\bm{\lambda}+\bm{\delta})
  +\mathcal{E}_1(\bm{\lambda}),
  \label{shapeinvJ}
\end{equation}
which is the (similarity transformed) shape invariance condition.
The constants $f^{\text{J}}_n$ and $b^{\text{J}}_n$
satisfy
\begin{equation}
  \mathcal{E}_n(\bm{\lambda})
  =f^{\text{J}}_n(\bm{\lambda})b^{\text{J}}_{n-1}(\bm{\lambda})
  \ \ (n\in\mathbb{Z}_{\geq 0}),
\end{equation}
and the energy eigenvalues $\mathcal{E}_n$ satisfy \eqref{Enshapeinv}.
Note that we have $f^{\text{J}}_n(\bm{\lambda})=\mathcal{E}_n(\bm{\lambda})$
and $b^{\text{J}}_n(\bm{\lambda})=1$ due to our normalization of $P_n(\eta)$.
Corresponding to the factorization \eqref{tHJ=BJFJ}, we have the following
relations \cite{os34}.
\begin{mythm}\cite{kls}
\label{thm:fbsrJ}
For the polynomials in \eqref{polyrdQMJ},
the forward and backward shift relations hold:
\begin{align}
  \mathcal{F}^{\text{\rm J}}(\bm{\lambda})P_n(\eta;\bm{\lambda})
  &=f^{\text{\rm J}}_n(\bm{\lambda})P_{n-1}(r \eta;\bm{\lambda}+\bm{\delta})
  \ \ (n\in\mathbb{Z}_{\geq 0}),
  \label{FJPn=}\\
  \mathcal{B}^{\text{\rm J}}(\bm{\lambda})
  P_{n-1}(r \eta;\bm{\lambda}+\bm{\delta})
  &=b^{\text{\rm J}}_{n-1}(\bm{\lambda})P_n(\eta;\bm{\lambda})
  \ \ (n\in\mathbb{Z}_{\geq 1}).
  \label{BJPn-1=}
\end{align}
\end{mythm}

\medskip

Let us comment on the Rodrigues-type formulas \cite{kls}.
The backward shift relations \eqref{BcPn-1=}, \eqref{BJPn-1=} and
$\check{P}_0(x)=P_0(\eta)=1$ give the following formulas,
\begin{align}
  \text{oQM,\,idQM}&:\ \check{P}_n(x;\bm{\lambda})
  =\prod_{j=0}^{n-1}b_{n-1-j}(\bm{\lambda}+j\bm{\delta})^{-1}\cdot
  \raprod{j=0}{n-1}\mathcal{B}(\bm{\lambda}+j\bm{\delta})\cdot 1,
  \label{Rodrigues1}\\
  \text{rdQM}&:\ \check{P}_n(x;\bm{\lambda})
  =\raprod{j=0}{n-1}\mathcal{B}(\bm{\lambda}+j\bm{\delta})\cdot 1,
  \label{Rodrigues1r}\\
  \text{rdQMJ}&:\ P_n(\eta;\bm{\lambda})=
  \raprod{j=0}{n-1}\mathcal{B}^{\text{J}}(\bm{\lambda}+j\bm{\delta})
  \Bigl|_{\eta\to r^j\eta}\cdot\,1,
  \label{Rodrigues1J}
\end{align}
where $\prod\limits_{k=1}^{\stackrel{n}{\longrightarrow}}a_k\eqdef
a_1a_2\cdots a_n$.
For oQM, idQM and rdQM, the Hamiltonian $\mathcal{H}(\bm{\lambda})$ is
factorized, $\mathcal{H}(\bm{\lambda})=\mathcal{A}(\bm{\lambda})^{\dagger}
\mathcal{A}(\bm{\lambda})$, and the operators $\mathcal{A}(\bm{\lambda})$ and
$\mathcal{A}(\bm{\lambda})^{\dagger}$ are related to the operators
$\mathcal{F}(\bm{\lambda})$ and $\mathcal{B}(\bm{\lambda})$ as follows
\cite{os24}:
\begin{align}
  \mathcal{F}(\bm{\lambda})&=\phi_0(x;\bm{\lambda}+\bm{\delta})^{-1}\circ
  \mathcal{A}(\bm{\lambda})\circ\phi_0(x;\bm{\lambda})
  \times\left\{
  \begin{array}{ll}
  1&:\text{oQM,\,idQM}\\
  \sqrt{B(0;\bm{\lambda})}&:\text{rdQM}
  \end{array}\right.,
  \label{FfromA}\\
  \mathcal{B}(\bm{\lambda})&=\phi_0(x;\bm{\lambda})^{-1}\circ
  \mathcal{A}(\bm{\lambda})^{\dagger}\circ\phi_0(x;\bm{\lambda}+\bm{\delta})
  \times\left\{
  \begin{array}{ll}
  1&:\text{oQM,\,idQM}\\
  \frac{1}{\sqrt{B(0;\bm{\lambda})}}&:\text{rdQM}
  \end{array}\right.,
  \label{BfromAd}
\end{align}
where $\phi_0(x;\bm{\lambda})$ is the ground state wavefunction, and
$\phi_0(x;\bm{\lambda})^2$ is the weight function of the orthogonal polynomials
$\check{P}_n(x;\bm{\lambda})$.
Explicit forms of $\phi_0(x;\bm{\lambda})$ are found in \cite{os13} for idQM,
\cite{os12} for rdQM, and $\phi_0(x;\bm{\lambda})=e^{w(x;\bm{\lambda})}$ for oQM.
By using \eqref{FfromA}--\eqref{BfromAd} and $\cF\phi_0(x;\bm{\lambda})=
\frac{d\eta(x)}{dx}\phi_0(x;\bm{\lambda}-\bm{\delta})$ for oQM,
the backward shift operators \eqref{BoQM}, \eqref{BidQM} and \eqref{BrdQM} are
rewritten as follows \cite{s16,os34}:
\begin{align}
  \text{oQM}&:\ \mathcal{B}(\bm{\lambda})=\phi_0(x;\bm{\lambda})^{-2}\circ
  \Bigl(-\cF\frac{d}{dx}\circ\Bigl(\frac{d\eta(x)}{dx}\Bigr)^{-1}\Bigr)\circ
  \phi_0(x;\bm{\lambda}+\bm{\delta})^2\n
  &\phantom{:\ \mathcal{B}(\bm{\lambda})\ }
  =-\cF\bigl(\phi_0(x;\bm{\lambda}-\bm{\delta})\phi_0(x;\bm{\lambda})\bigr)^{-1}
  \circ\frac{d}{d\eta}\circ
  \bigl(\phi_0(x;\bm{\lambda})\phi_0(x;\bm{\lambda}+\bm{\delta})\bigr),\\
  \text{idQM}&:\ \mathcal{B}(\bm{\lambda})=\phi_0(x;\bm{\lambda})^{-2}\circ
  \mathcal{D}\circ\phi_0(x;\bm{\lambda}+\bm{\delta})^2,\quad
  \mathcal{D}\eqdef-i(e^{\frac{\gamma}{2}p}-e^{-\frac{\gamma}{2}p})
  \varphi(x)^{-1},
  \label{calDidQM}\\
  \text{rdQM}&:\ \mathcal{B}(\bm{\lambda})=\phi_0(x;\bm{\lambda})^{-2}\circ
  \mathcal{D}(\bm{\lambda})\circ\phi_0(x;\bm{\lambda}+\bm{\delta})^2,\quad
  \mathcal{D}(\bm{\lambda})\eqdef(1-e^{-\partial})\varphi(x;\bm{\lambda})^{-1}.
\end{align}
Although the function $\varphi(x)$ depends on the parameters $\bm{\lambda}$
for c$q$H and $q$MP cases in idQM, we suppress writing $\bm{\lambda}$ dependence
in $\mathcal{D}$ \eqref{calDidQM}, because $\varphi(x;\bm{\lambda})$ is
invariant under the shift $\bm{\lambda}\to\bm{\lambda}+\bm{\delta}$.
Then \eqref{Rodrigues1} and \eqref{Rodrigues1r} become
\begin{align}
  \text{oQM}&:\ \check{P}_n(x;\bm{\lambda})=(-\cF)^n
  \prod_{j=0}^{n-1}b_{n-1-j}(\bm{\lambda}+j\bm{\delta})^{-1}
  \label{Rodrigues2o}\\
  &\phantom{:\ \check{P}_n(x;\bm{\lambda})=}\times
  \bigl(\phi_0(x;\bm{\lambda}-\bm{\delta})\phi_0(x;\bm{\lambda})\bigr)^{-1}
  \Bigl(\frac{d}{d\eta}\Bigr)^n\cdot
  \phi_0\bigl(x;\bm{\lambda}+(n-1)\bm{\delta}\bigr)
  \phi_0(x;\bm{\lambda}+n\bm{\delta}),\n
  \text{idQM}&:\ \check{P}_n(x;\bm{\lambda})=
  \prod_{j=0}^{n-1}b_{n-1-j}(\bm{\lambda}+j\bm{\delta})^{-1}\cdot
  \phi_0(x;\bm{\lambda})^{-2}\,\mathcal{D}^n\cdot
  \phi_0(x;\bm{\lambda}+n\bm{\delta})^2,
  \label{Rodrigues2i}\\
  \text{rdQM}&:\ \check{P}_n(x;\bm{\lambda})=\phi_0(x;\bm{\lambda})^{-2}
  \raprod{j=0}{n-1}\mathcal{D}(\bm{\lambda}+j\bm{\delta})\cdot
  \phi_0(x;\bm{\lambda}+n\bm{\delta})^2.
  \label{Rodrigues2r}
\end{align}
For rdQM, there are five sinusoidal coordinates \cite{os12}:
(\romannumeral1) $\eta(x)=x$,
(\romannumeral2) $\eta(x;\bm{\lambda})=x(x+d)$,
(\romannumeral3) $\eta(x)=1-q^x$,
(\romannumeral4) $\eta(x)=q^{-x}-1$,
(\romannumeral5) $\eta(x;\bm{\lambda})=(q^{-x}-1)(1-dq^x)$.
For (\romannumeral1), (\romannumeral3) and (\romannumeral4),
the operator $\mathcal{D}(\bm{\lambda})$ does not depend on $\bm{\lambda}$,
and we have $\prod\limits_{j=0}^{\stackrel{n-1}{\longrightarrow}}
\mathcal{D}(\bm{\lambda}+j\bm{\delta})=\mathcal{D}^n$.
We remark that \eqref{Rodrigues2r} can be rewritten as
\begin{equation}
  \check{P}_n(x;\bm{\lambda})=\Bigl(
  \frac{\phi_0(x;\bm{\lambda})^2}{\varphi(x;\bm{\lambda}-\bm{\delta})}
  \Bigr)^{-1}
  \raprod{j=0}{n-1}\varphi\bigl(x;\bm{\lambda}+(j-1)\bm{\delta}\bigr)^{-1}
  (1-e^{-\partial})\cdot
  \frac{\phi_0(x;\bm{\lambda}+n\bm{\delta})^2}
  {\varphi(x;\bm{\lambda}+(n-1)\bm{\delta})},
\end{equation}
and the expressions for the formulas (9.2.10), (14.2.11) etc.\ in
\cite{kls} are somewhat ambiguous (about the shift of parameters).
For rdQMJ, direct calculation shows that the backward shift operators
\eqref{BJrdQMJ} are rewritten as follows (cf.\,\cite{os34}):
\begin{align}
  \mathcal{B}(\bm{\lambda})=A'(\bm{\lambda})
  \phi^{\text{J}}_0(\eta;\bm{\lambda})^{-2}\circ\mathcal{D}^{\text{J}}\circ
  \phi^{\text{J}}_0(r\eta;\bm{\lambda}+\bm{\delta})^2,\quad
  \mathcal{D}^{\text{J}}\eqdef(1-q^{-\eta\frac{d}{d\eta}})\eta^{-1},
  \label{calDrdQMJ}
\end{align}
where the function $\phi^{\text{J}}_0(\eta;\bm{\lambda})$ and the constant
$A'(\bm{\lambda})$ are given by
\begin{align}
  &\text{b$q$J,\,b$q$L,\,ACS$\I$,\,d$q$He$\I$,\,d$q$He$\II$}
  :\ \ \phi^{\text{J}}_0(\eta;\bm{\lambda})^2\eqdef
  \eta\prod_{j=0}^{\infty}\frac{qD^{\text{J}}(q^{j+1}\eta;\bm{\lambda})}
  {B^{\text{J}}(q^j\eta;\bm{\lambda})},\\
  &\text{$q$L}:\ \ \phi^{\text{J}}_0(\eta;\bm{\lambda})^2
  \eqdef\frac{\eta^{\lambda_1+1}}{c^{\lambda_1+1}}\frac{1}{(-\eta;q)_{\infty}},
  \ \ \text{namely}\ \phi^{\text{J}}_0(cq^x;\bm{\lambda})^2
  =\frac{(aq)^x}{(-cq^x;q)_{\infty}},\\
  &\text{SW}:\ \ \phi^{\text{J}}_0(\eta;\bm{\lambda})^2
  \eqdef(c^{-1}\eta)^{\frac12(\frac{\log\eta}{\log q}+\lambda_1+1)},
  \ \ \text{namely}\ \phi^{\text{J}}_0(cq^x;\bm{\lambda})^2
  =c^xq^{\frac12x(x+1)},\\
  &A'(\bm{\lambda})=A(\bm{\lambda})^{-1}\times\left\{
  \begin{array}{ll}
  (-ac)&:\text{b$q$J}\\
  (-ab)&:\text{b$q$L}\\
  (-aq^{-1})&:\text{ASC$\I$}\\
  q^{-1}&:\text{d$q$He$\I$,\,d$q$He$\II$}\\
  c&:\text{$q$L}\\
  c^2&:\text{SW}
  \end{array}\right..
\end{align}
Here the parameters $\bm{\lambda}$ are extended to
$q^{\bm{\lambda}}=(a,c)$ with $\bm{\delta}=(1,1)$ for $q$L and
$q^{\bm{\lambda}}=c$ with $\bm{\delta}=2$ for SW.
Then \eqref{Rodrigues1J} becomes
\begin{equation}
  P_n(\eta;\bm{\lambda})=r^{-\binom{n}{2}}\prod_{j=0}^{n-1}
  A'(\bm{\lambda}+j\bm{\delta})\cdot
  \phi^{\text{J}}_0(\eta;\bm{\lambda})^{-2}
  \bigl(\mathcal{D}^{\text{J}}\bigr)^n\cdot
  \phi^{\text{J}}_0(r^n\eta;\bm{\lambda}+n\bm{\delta})^2.
\end{equation}
We remark that
\begin{equation}
  \bigl((1-q^{-\eta\frac{d}{d\eta}})\eta^{-1}\bigr)^n
  \phi^{\text{J}}_0(q^n\eta;\bm{\lambda}+n\bm{\delta})^2
  =(-1)^nq^{\frac12n(n+1)}
  \bigl((1-q^{\eta\frac{d}{d\eta}})\eta^{-1}\bigr)^n
  \phi^{\text{J}}_0(\eta;\bm{\lambda}+n\bm{\delta})^2
\end{equation}
for b$q$J,\,b$q$L,\,d$q$He$\II$ and $qL$, and
$\phi^{\text{J}}_0(q^{2n}\eta;\bm{\lambda}+n\bm{\delta})^2
=(c^{-1}\eta)^{2n}\phi^{\text{J}}_0(\eta;\bm{\lambda})^2$ for SW.

%%%%%%%%%%%%%%%%%%%%%%%%%%%%%%%%%%%%%%%%%%%%%%%%%%%%%%%%%%%%%%%
%                                                             %
%  4. Another type of Forward and Backward Shift Relations    %
%                                                             %
%%%%%%%%%%%%%%%%%%%%%%%%%%%%%%%%%%%%%%%%%%%%%%%%%%%%%%%%%%%%%%%
\section{Another type of Forward and Backward Shift Relations}
\label{sec:newfbsr}

In this section, based on other factorizations of $\widetilde{\mathcal{H}}$
and $\widetilde{\mathcal{H}}^{\text{J}}$, we present another type of forward
and backward shift relations.

%%%%%%%%%%%%%%%%%%%%%%%%%%%%%%%%%%%%%%%%%%%%%%%%%%%%%
%                                                   %
% 4.1 Polynomials in oQM systems                    %
%                                                   %
%%%%%%%%%%%%%%%%%%%%%%%%%%%%%%%%%%%%%%%%%%%%%%%%%%%%%
\subsection{Polynomials in oQM systems}
\label{sec:newoQM}
\setcounter{mybangouI}{3}
\setcounter{mybangouII}{1}
\setcounter{mythm}{0}

For the oQM systems described by the polynomials \eqref{polyoQM}
(except He and B),
let us define the operators $\tilde{\mathcal{F}}(\bm{\lambda})$ and
$\tilde{\mathcal{B}}(\bm{\lambda})$.
For J case, they are given by
\begin{align}
%%%%%%%%%%%%%%%%%%%%%%%%%%%%%%%%
% Jacobi (J)                   %
%%%%%%%%%%%%%%%%%%%%%%%%%%%%%%%%
  \text{(a)}&:
  \ \tilde{\mathcal{F}}(\bm{\lambda})\eqdef\frac12\tan x\frac{d}{dx}+g-\frac12
  \ \Bigl(=-(1-\eta)\frac{d}{d\eta}+g-\frac12\Bigr),
  \label{tFoQM:J(a)}\\
  &\ \ \ \ \tilde{\mathcal{B}}(\bm{\lambda})\eqdef
  -\frac12\cot x\frac{d}{dx}+h+\frac12
  \ \Bigl(=(1+\eta)\frac{d}{d\eta}+h+\frac12\Bigr),
  \label{tBoQM:J(a)}\\
  \text{(b)}&:
  \ \tilde{\mathcal{F}}(\bm{\lambda})\eqdef-\frac12\cot x\frac{d}{dx}+h-\frac12
  \ \Bigl(=(1+\eta)\frac{d}{d\eta}+h-\frac12\Bigr),
  \label{tFoQM:J(b)}\\
  &\ \ \ \ \tilde{\mathcal{B}}(\bm{\lambda})\eqdef
  \frac12\tan x\frac{d}{dx}+g+\frac12
  \ \Bigl(=-(1-\eta)\frac{d}{d\eta}+g+\frac12\Bigr),
  \label{tBoQM:J(b)}
\end{align}
and the constants $\tilde{f}_n(\bm{\lambda})$, $\tilde{b}_n(\bm{\lambda})$
and $\bm{\bar{\delta}}$ are given by
\begin{align}
%%%%%%%%%%%%%%%%%%%%%%%%%%%%%%%%
% Jacobi (J)                   %
%%%%%%%%%%%%%%%%%%%%%%%%%%%%%%%%
  \text{(a)}&:
  \ \tilde{f}_n(\bm{\lambda})=n+g-\tfrac12,
  \ \ \tilde{b}_n(\bm{\lambda})=n+h+\tfrac12,
  \ \ \bm{\bar{\delta}}=(1,-1),\\
  \text{(b)}&:
  \ \tilde{f}_n(\bm{\lambda})=n+h-\tfrac12,
  \ \ \tilde{b}_n(\bm{\lambda})=n+g+\tfrac12,
  \ \ \bm{\bar{\delta}}=(-1,1).
\end{align}
For L and pJ cases, see Appendix \ref{app:oQM}.

Then we can show that
\begin{align}
  &\widetilde{\mathcal{H}}(\bm{\lambda})
  =4\bigl(\tilde{\mathcal{B}}(\bm{\lambda})\tilde{\mathcal{F}}(\bm{\lambda})
  -\tilde{f}_0(\bm{\lambda})\tilde{b}_0(\bm{\lambda})\bigr),
  \label{tH=tBtF..oQM}\\
  &\mathcal{E}_n(\bm{\lambda})
  =4\bigl(\tilde{f}_n(\bm{\lambda})\tilde{b}_n(\bm{\lambda})
  -\tilde{f}_0(\bm{\lambda})\tilde{b}_0(\bm{\lambda})\bigr)
  \ \ (n\in\mathbb{Z}_{\geq 0}).
\end{align}
Corresponding to this factorization \eqref{tH=tBtF..oQM},
the following relations are obtained by direct calculation.
\begin{mythm}
\label{thm:newfbsr_oQM}
For the polynomials in \eqref{polyoQM} (except {\rm He} and {\rm B}),
the following forward and backward shift relations hold for
$n\in\mathbb{Z}_{\geq 0}$,
\begin{align}
  \tilde{\mathcal{F}}(\bm{\lambda})\check{P}_n(x;\bm{\lambda})
  &=\tilde{f}_n(\bm{\lambda})\check{P}_n(x;\bm{\lambda}-\bm{\bar{\delta}}),
  \label{tFcPn=oQM}\\
  \tilde{\mathcal{B}}(\bm{\lambda})\check{P}_n(x;\bm{\lambda}-\bm{\bar{\delta}})
  &=\tilde{b}_n(\bm{\lambda})\check{P}_n(x;\bm{\lambda}).
  \label{tBcPn=oQM}
\end{align}
\end{mythm}

\remark
We think that these identities \eqref{tFcPn=oQM}--\eqref{tBcPn=oQM} may be
known formulas but this interpretation is new.
For example, \eqref{tFcPn=oQM} for \eqref{tFoQM:J(a)} can be obtained
by differentiating the integral formula (1.18) with $\mu=1$ in \cite{af69}.

\remark
Two formulas with $\bm{\bar{\delta}}$ and $-\bm{\bar{\delta}}$ are equivalent
by interchanging $\tilde{\mathcal{F}}$ and $\tilde{\mathcal{B}}$, e.g.
\eqref{tFcPn=oQM} and \eqref{tBcPn=oQM} for L\,(b) agree with
\eqref{tBcPn=oQM} and \eqref{tFcPn=oQM} for L\,(a) with the replacement
$g\to g+1$, respectively.
For He and B, we do not have new factorization \eqref{tH=tBtF..oQM} and
another type of forward and backward shift relations
\eqref{tFcPn=oQM}--\eqref{tBcPn=oQM}.

%%%%%%%%%%%%%%%%%%%%%%%%%%%%%%%%%%%%%%%%%%%%%%%%%%%%%
%                                                   %
% 4.2 Polynomials in idQM systems                   %
%                                                   %
%%%%%%%%%%%%%%%%%%%%%%%%%%%%%%%%%%%%%%%%%%%%%%%%%%%%%
\subsection{Polynomials in idQM systems}
\label{sec:newidQM}
\addtocounter{mybangouII}{1}
\setcounter{myrembangou}{0}

For the idQM systems described by the polynomials \eqref{polyidQM},
let us define the operators $\tilde{\mathcal{F}}(\bm{\lambda})$ and
$\tilde{\mathcal{B}}(\bm{\lambda})$ as follows:
\begin{align}
  \tilde{\mathcal{F}}(\bm{\lambda})&\eqdef
  V_1(x+i\tfrac{\gamma}{2};\bm{\lambda})e^{\frac{\gamma}{2}p}
  +V_1^*(x-i\tfrac{\gamma}{2};\bm{\lambda})e^{-\frac{\gamma}{2}p},
  \label{tFidQM}\\
  \tilde{\mathcal{B}}(\bm{\lambda})&\eqdef
  V_2(x;\bm{\lambda})e^{\frac{\gamma}{2}p}
  +V_2^*(x;\bm{\lambda})e^{-\frac{\gamma}{2}p},
  \label{tBidQM}
\end{align}
where the potential functions $V_1(x;\bm{\lambda})$ and $V_2(x;\bm{\lambda})$
satisfy
\begin{equation}
  V(x;\bm{\lambda})=V_1(x;\bm{\lambda})V_2(x;\bm{\lambda}).
  \label{V=V1V2}
\end{equation}
For AW case, their explicit forms are given by
\begin{align}
%%%%%%%%%%%%%%%%%%%%%%%%%%%%%%%%
% Askey-Wilson (AW)            %
%%%%%%%%%%%%%%%%%%%%%%%%%%%%%%%%
  &\text{Assume $\{a_j^*,a_k^*\}=\{a_j,a_k\}$ (as a set) and
  set $\{l,m\}=\{1,2,3,4\}\backslash\{j,k\}$,}\n
  &V_i(x;\bm{\lambda})=V^{(j,k)}_i(x;\bm{\lambda})\ \ (i=1,2),\n
  &V_1(x;\bm{\lambda})=\frac{(1-a_je^{ix})(1-a_ke^{ix})}{1-qe^{2ix}},
  \ \ V_2(x;\bm{\lambda})=\frac{(1-a_le^{ix})(1-a_me^{ix})}{1-e^{2ix}},
  \label{V1V2:AW}
\end{align}
and the constants $\tilde{f}_n(\bm{\lambda})$, $\tilde{b}_n(\bm{\lambda})$
and $\bm{\bar{\delta}}$ as given by
\begin{align}
%%%%%%%%%%%%%%%%%%%%%%%%%%%%%%%%
% Askey-Wilson (AW)            %
%%%%%%%%%%%%%%%%%%%%%%%%%%%%%%%%
  &\tilde{f}_n(\bm{\lambda})=q^{-\frac{n}{2}}(1-a_ja_kq^{n-1}),
  \ \ \tilde{b}_n(\bm{\lambda})=q^{-\frac{n}{2}}(1-a_la_mq^n),\n
  &(\bm{\bar{\delta}})_j=(\bm{\bar{\delta}})_k=\tfrac12,
  \ \ (\bm{\bar{\delta}})_l=(\bm{\bar{\delta}})_m=-\tfrac12.
\end{align}
For other cases, see Appendix \ref{app:idQM}.

Then we can show that $V_1(x)$ and $V_2(x)$ satisfy
\begin{equation}
  V_1(x+i\gamma;\bm{\lambda})V_2^*(x;\bm{\lambda})
  +V_1^*(x-i\gamma;\bm{\lambda})V_2(x;\bm{\lambda})
  -\tilde{f}_0(\bm{\lambda})\tilde{b}_0(\bm{\lambda})
  =-V(x;\bm{\lambda})-V^*(x;\bm{\lambda}),
  \label{V1V*2...}
\end{equation}
and the constants $\tilde{f}_n$ and $\tilde{b}_n$ satisfy
\begin{equation}
  \mathcal{E}_n(\bm{\lambda})
  =\tilde{f}_n(\bm{\lambda})\tilde{b}_n(\bm{\lambda})
  -\tilde{f}_0(\bm{\lambda})\tilde{b}_0(\bm{\lambda})
  \ \ (n\in\mathbb{Z}_{\geq 0}).
  \label{En=ftnbtn-..}
\end{equation}
The relation \eqref{V1V*2...} gives other factorizations of
$\widetilde{\mathcal{H}}(\bm{\lambda})$ \eqref{tHidQM},
\begin{equation}
  \widetilde{\mathcal{H}}(\bm{\lambda})
  =\tilde{\mathcal{B}}(\bm{\lambda})\tilde{\mathcal{F}}(\bm{\lambda})
  -\tilde{f}_0(\bm{\lambda})\tilde{b}_0(\bm{\lambda}).
  \label{tH=tBtF..idQM}
\end{equation}
Corresponding to this factorization \eqref{tH=tBtF..idQM},
we obtain the following relations.
\begin{mythm}
\label{thm:newfbsr_idQM}
For the polynomials in \eqref{polyidQM},
the following forward and backward shift relations hold for
$n\in\mathbb{Z}_{\geq 0}$,
\begin{align}
  \tilde{\mathcal{F}}(\bm{\lambda})\check{P}_n(x;\bm{\lambda})
  &=\tilde{f}_n(\bm{\lambda})\check{P}_n(x;\bm{\lambda}-\bm{\bar{\delta}}),
  \label{tFcPn=idQM}\\
  \tilde{\mathcal{B}}(\bm{\lambda})\check{P}_n(x;\bm{\lambda}-\bm{\bar{\delta}})
  &=\tilde{b}_n(\bm{\lambda})\check{P}_n(x;\bm{\lambda}).
  \label{tBcPn=idQM}
\end{align}
\end{mythm}
Proof:
It is sufficient to show \eqref{tFcPn=idQM}, because \eqref{tHcPn=} and
\eqref{En=ftnbtn-..}--\eqref{tFcPn=idQM} imply \eqref{tBcPn=idQM}.
Taking AW \eqref{V1V2:AW} with $(j,k)=(1,2)$ as an example, let us prove
\eqref{tFcPn=idQM}.
It is shown by direct calculation:
\begin{align*}
  &\quad\tilde{\mathcal{F}}(\bm{\lambda})\check{P}_n(x;\bm{\lambda})\\
  &=\frac{(1-a_1q^{-\frac12}e^{ix})(1-a_2q^{-\frac12}e^{ix})}{1-e^{2ix}}
  \frac{(a_1a_2,a_1a_3,a_1a_4\,;q)_n}{a_1^n}\\
  &\qquad\times
  {}_4\phi_3\Bigl(
  \genfrac{}{}{0pt}{}{q^{-n},\,a_1a_2a_3a_4q^{n-1},\,a_1q^{\frac12}e^{ix},
  \,a_1q^{-\frac12}e^{-ix}}
  {a_1a_2,\,a_1a_3,\,a_1a_4}\Bigm|q\,;q\Bigr)\\
  &\quad+\frac{(1-a_1q^{-\frac12}e^{-ix})(1-a_2q^{-\frac12}e^{-ix})}{1-e^{-2ix}}
  \frac{(a_1a_2,a_1a_3,a_1a_4\,;q)_n}{a_1^n}\\
  &\qquad\times
  {}_4\phi_3\Bigl(
  \genfrac{}{}{0pt}{}{q^{-n},\,a_1a_2a_3a_4q^{n-1},\,a_1q^{-\frac12}e^{ix},
  \,a_1q^{\frac12}e^{-ix}}
  {a_1a_2,\,a_1a_3,\,a_1a_4}\Bigm|q\,;q\Bigr)\\
  &=\frac{(a_1a_2,a_1a_3,a_1a_4\,;q)_n}{a_1^n(1-e^{2ix})}
  \sum_{k=0}^n\frac{(q^{-n},a_1a_2a_3a_4q^{n-1},a_1q^{-\frac12}e^{ix},
  a_1q^{-\frac12}e^{-ix}\,;q)_k}{(a_1a_2,a_1a_3,a_1a_4\,;q)_k}
  \frac{q^k}{(q\,;q)_k}\\
  &\qquad\times\bigl((1-a_1e^{ix}q^{k-\frac12})(1-a_2q^{-\frac12}e^{ix})
  -e^{2ix}(1-a_1e^{-ix}q^{k-\frac12})(1-a_2q^{-\frac12}e^{-ix})\bigr)\\
  &=\frac{(a_1a_2,a_1a_3,a_1a_4\,;q)_n}{a_1^n(1-e^{2ix})}
  \sum_{k=0}^n\frac{(q^{-n},a_1a_2a_3a_4q^{n-1},a_1q^{-\frac12}e^{ix},
  a_1q^{-\frac12}e^{-ix}\,;q)_k}{(a_1a_2,a_1a_3,a_1a_4\,;q)_k}
  \frac{q^k}{(q\,;q)_k}\\
  &\qquad\times(1-a_1a_2q^{k-1})(1-e^{2ix})\\
  &=\frac{(a_1a_2,a_1a_3,a_1a_4\,;q)_n}{a_1^n}
  \sum_{k=0}^n(1-a_1a_2q^{-1})
  \frac{(q^{-n},a_1a_2a_3a_4q^{n-1},a_1q^{-\frac12}e^{ix},
  a_1q^{-\frac12}e^{-ix}\,;q)_k}{(a_1a_2q^{-1},a_1a_3,a_1a_4\,;q)_k}
  \frac{q^k}{(q\,;q)_k}\\
  &=q^{-\frac{n}{2}}(1-a_1a_2q^{n-1})\\
  &\qquad\times
  \frac{(a_1a_2q^{-1},a_1a_3,a_1a_4\,;q)_n}{(a_1q^{-\frac12})^n}
  \sum_{k=0}^n\frac{(q^{-n},a_1a_2a_3a_4q^{n-1},a_1q^{-\frac12}e^{ix},
  a_1q^{-\frac12}e^{-ix}\,;q)_k}{(a_1a_2q^{-1},a_1a_3,a_1a_4\,;q)_k}
  \frac{q^k}{(q\,;q)_k}\\
  &=\tilde{f}_n(\bm{\lambda})\check{P}_n(x;\bm{\lambda}-\bm{\bar{\delta}}).
\end{align*}
The other cases are proved in the same way.
\hfill\fbox{}

\medskip

\remark
After submitting this manuscript to arXiv (arXiv:2301.00678),
when I gave a talk at ICIAM 2023 Tokyo (August 2023),
Satoshi Tsujimoto informed me of the paper by Kalnins and Miller \cite{km89},
in which the formulas \eqref{tFcPn=idQM}--\eqref{tBcPn=idQM} for Askey-Wilson
case \eqref{V1V2:AW} with $(j,k)=(1,2)$ are given from a different point of
view.
Their operators $\tau=\tau^{(a,b,c,d)}$,
$\tau^*=\tau^{*(aq^{\frac12},bq^{\frac12},cq^{\frac12},dq^{\frac12})}$,
$\mu=\mu^{(a,b,c,d)}$ and
$\mu^*=\mu^{(cq^{\frac12},dq^{\frac12},aq^{-\frac12},bq^{-\frac12})}$
correspond to our
$\mathcal{F}(\bm{\lambda})$, $\mathcal{B}(\bm{\lambda})$,
$\tilde{\mathcal{F}}(\bm{\lambda})$ and $\tilde{\mathcal{B}}(\bm{\lambda})$
with $(a_1,a_2,a_3,a_4)=(a,b,c,d)$, respectively.

\remark
Two formulas with $\bm{\bar{\delta}}$ and $-\bm{\bar{\delta}}$ are equivalent
by interchanging $\tilde{\mathcal{F}}$ and $\tilde{\mathcal{B}}$, e.g.
\eqref{tFcPn=idQM} and \eqref{tBcPn=idQM} for cH\,(b) agree with
\eqref{tBcPn=idQM} and \eqref{tFcPn=idQM} for cH\,(a) with the replacements
$a_1\to a_1+\frac12$ and $a_2\to a_2-\frac12$, respectively.

\remark\label{*inv}\!\!\!
The relations \eqref{tFcPn=idQM}--\eqref{tBcPn=idQM} are invariant under the
$*$-operation.
In contrast to the $x$-shift relations studied in \cite{os40} (see
Theorem\,\ref{thm:newfbsr_rdQM}), the coordinate $x$ is not shifted, and only
the parameters $\bm{\lambda}$ are shifted.
We choose the operators $\tilde{\mathcal{F}}(\bm{\lambda})$ and
$\tilde{\mathcal{B}}(\bm{\lambda})$ \eqref{tFidQM}--\eqref{tBidQM} to respect
this $*$-operation invariance.
See also Remark\,\ref{AWqR}.

\remark\label{simodoki:idQM}\!\!\!
We can show that
\begin{align}
  &V_1(x+i\tfrac{\gamma}{2};\bm{\lambda})V_2(x-i\tfrac{\gamma}{2};\bm{\lambda})
  =V_1(x;\bm{\lambda}-\bm{\bar{\delta}})
  V_2(x;\bm{\lambda}-\bm{\bar{\delta}}),\n
  &V_1(x+i\tfrac{\gamma}{2};\bm{\lambda})V_2^*(x-i\tfrac{\gamma}{2};\bm{\lambda})
  +V_1^*(x-i\tfrac{\gamma}{2};\bm{\lambda})V_2(x+i\tfrac{\gamma}{2};\bm{\lambda})
  -\tilde{f}_0(\bm{\lambda})\tilde{b}_0(\bm{\lambda})\\
  &=V_1(x+i\gamma;\bm{\lambda}-\bm{\bar{\delta}})
  V_2^*(x;\bm{\lambda}-\bm{\bar{\delta}})
  +V_1^*(x-i\gamma;\bm{\lambda}-\bm{\bar{\delta}})
  V_2(x;\bm{\lambda}-\bm{\bar{\delta}})
  -\tilde{f}_0(\bm{\lambda}-\bm{\bar{\delta}})
  \tilde{b}_0(\bm{\lambda}-\bm{\bar{\delta}}),
  \nonumber
\end{align}
which imply
\begin{equation}
  \tilde{\mathcal{F}}(\bm{\lambda})\tilde{\mathcal{B}}(\bm{\lambda})
  -\tilde{f}_0(\bm{\lambda})\tilde{b}_0(\bm{\lambda})
  =\tilde{\mathcal{B}}(\bm{\lambda}-\bm{\bar{\delta}})
  \tilde{\mathcal{F}}(\bm{\lambda}-\bm{\bar{\delta}})
  -\tilde{f}_0(\bm{\lambda}-\bm{\bar{\delta}})
  \tilde{b}_0(\bm{\lambda}-\bm{\bar{\delta}}).
  \label{tFtB-..:idQM}
\end{equation}

%%%%%%%%%%%%%%%%%%%%%%%%%%%%%%%%%%%%%%%%%%%%%%%%%%%%%
%                                                   %
% 4.3 Polynomials in rdQM systems                   %
%                                                   %
%%%%%%%%%%%%%%%%%%%%%%%%%%%%%%%%%%%%%%%%%%%%%%%%%%%%%
\subsection{Polynomials in rdQM systems}
\label{sec:newrdQM}
\addtocounter{mybangouII}{1}
\setcounter{myrembangou}{0}

For the rdQM systems described by the polynomials
\eqref{polyrdQMf}--\eqref{polyrdQMsi}
(except C and $q$B),
let us define the operators $\tilde{\mathcal{F}}(\bm{\lambda})$ and
$\tilde{\mathcal{B}}(\bm{\lambda})$
as follows:
\begin{align}
  \tilde{\mathcal{F}}(\bm{\lambda})&\eqdef
  D_1(x+1;\bm{\lambda})+B_1(x;\bm{\lambda})e^{\partial},
  \label{tFrdQM}\\
  \tilde{\mathcal{B}}(\bm{\lambda})&\eqdef
  B_2(x;\bm{\lambda})+D_2(x;\bm{\lambda})e^{-\partial},
  \label{tBrdQM}
\end{align}
where the potential functions $B_1(x;\bm{\lambda}),B_2(x;\bm{\lambda}),
D_1(x;\bm{\lambda})$ and $D_2(x;\bm{\lambda})$ satisfy
\begin{equation}
  B(x;\bm{\lambda})=B_1(x;\bm{\lambda})B_2(x;\bm{\lambda}),\quad
  D(x;\bm{\lambda})=D_1(x;\bm{\lambda})D_2(x;\bm{\lambda}).
  \label{B=B1B2}
\end{equation}
For $q$R case, their explicit forms are given by
\begin{align}
%%%%%%%%%%%%%%%%%%%%%%%%%%%%%%%%
% q-Racah (qR)                 %
%%%%%%%%%%%%%%%%%%%%%%%%%%%%%%%%
  \text{(a)}&:
  \ B_1(x;\bm{\lambda})=-\frac{(1-q^{x-N})(1-dq^x)}
  {(q^{-N-1}-1)(1-dq^{2x+1})},
  \ \ B_2(x;\bm{\lambda})=\frac{(q^{-N-1}-1)(1-bq^x)(1-cq^x)}{1-dq^{2x}},\n
  &\ \ \,D_1(x;\bm{\lambda})=\frac{(1-dq^{x+N})(1-q^x)}
  {(1-q^{N+1})(1-dq^{2x-1})},
  \ D_2(x;\bm{\lambda})=-\tilde{d}\,
  \frac{(1-q^{N+1})(1-b^{-1}dq^x)(1-c^{-1}dq^x)}{1-dq^{2x}},
  \label{B1B2D1D2:qR}\\
  \text{(b)}&:
  \ B_1(x;\bm{\lambda})=\frac{(1-bq^x)(1-dq^x)}{(1-bq^{-1})(1-dq^{2x+1})},
  \ \ B_2(x;\bm{\lambda})=-\frac{(1-bq^{-1})(1-q^{x-N})(1-cq^x)}{1-dq^{2x}},\n
  &\ \ \ D_1(x;\bm{\lambda})=\frac{(1-b^{-1}dq^x)(1-q^x)}
  {(1-b^{-1}q)(1-dq^{2x-1})},
  \ D_2(x;\bm{\lambda})=-\tilde{d}\,
  \frac{(1-b^{-1}q)(1-dq^{x+N})(1-c^{-1}dq^x)}{1-dq^{2x}},\\
  \text{(c)}&:
  \ B_1(x;\bm{\lambda})=\frac{(1-cq^x)(1-dq^x)}{(1-cq^{-1})(1-dq^{2x+1})},
  \ \ B_2(x;\bm{\lambda})=-\frac{(1-cq^{-1})(1-q^{x-N})(1-bq^x)}{1-dq^{2x}},\n
  &\ \ \ D_1(x;\bm{\lambda})=\frac{(1-c^{-1}dq^x)(1-q^x)}
  {(1-c^{-1}q)(1-dq^{2x-1})},
  \ D_2(x;\bm{\lambda})=-\tilde{d}\,
  \frac{(1-c^{-1}q)(1-dq^{x+N})(1-b^{-1}dq^x)}{1-dq^{2x}},\\
  \text{(d)}&:
  \ B_1(x;\bm{\lambda})=-\frac{(1-c)(1-q^{x-N})(1-bq^x)}{1-dq^{2x+1}},
  \ \ B_2(x;\bm{\lambda})=\frac{(1-cq^x)(1-dq^x)}{(1-c)(1-dq^{2x})},\\
  &\ \ \ D_1(x;\bm{\lambda})=-\tilde{d}\,
  \frac{(1-c^{-1})(1-dq^{x+N})(1-b^{-1}dq^x)}{1-dq^{2x-1}},
  \ \ D_2(x;\bm{\lambda})=\frac{(1-c^{-1}dq^x)(1-q^x)}
  {(1-c^{-1})(1-dq^{2x})},\n
  \text{(e)}&:
  \ B_1(x;\bm{\lambda})=-\frac{(1-b)(1-q^{x-N})(1-cq^x)}{1-dq^{2x+1}},
  \ \ B_2(x;\bm{\lambda})=\frac{(1-bq^x)(1-dq^x)}{(1-b)(1-dq^{2x})},\\
  &\ \ \ D_1(x;\bm{\lambda})=-\tilde{d}\,
  \frac{(1-b^{-1})(1-dq^{x+N})(1-c^{-1}dq^x)}{1-dq^{2x-1}},
  \ \ D_2(x;\bm{\lambda})=\frac{(1-b^{-1}dq^x)(1-q^x)}
  {(1-b^{-1})(1-dq^{2x})},\n
  \text{(f)}&:
  \ B_1(x;\bm{\lambda})=\frac{(q^{-N}-1)(1-bq^x)(1-cq^x)}{1-dq^{2x+1}},
  \ \ B_2(x;\bm{\lambda})=-\frac{(1-q^{x-N})(1-dq^x)}{(q^{-N}-1)(1-dq^{2x})},\\
  &\ \ \ D_1(x;\bm{\lambda})=-\tilde{d}\,
  \frac{(1-q^N)(1-b^{-1}dq^x)(1-c^{-1}dq^x)}{1-dq^{2x-1}},
  \ \ D_2(x;\bm{\lambda})=\frac{(1-dq^{x+N})(1-q^x)}{(1-q^N)(1-dq^{2x})},
  \nonumber
\end{align}
and the constants $\tilde{f}_n(\bm{\lambda})$, $\tilde{b}_n(\bm{\lambda})$
and $\bm{\bar{\delta}}$ are given by
\begin{align}
%%%%%%%%%%%%%%%%%%%%%%%%%%%%%%%%
% q-Racah (qR)                 %
%%%%%%%%%%%%%%%%%%%%%%%%%%%%%%%%
  \text{(a)}&:
  \ \tilde{f}_n(\bm{\lambda})=1,
  \ \ \tilde{b}_n(\bm{\lambda})
  =\mathcal{E}_{N+1}(\bm{\lambda})-\mathcal{E}_n(\bm{\lambda}),
  \ \bm{\bar{\delta}}=(1,0,0,1),\\
  \text{(b)}&:
  \ \tilde{f}_n(\bm{\lambda})=1,
  \ \ \tilde{b}_n(\bm{\lambda})=-q^{-n}(1-bq^{n-1})(1-cd^{-1}q^{n-N}),
  \ \bm{\bar{\delta}}=(0,1,0,1),\!\!\\
  \text{(c)}&:
  \ \tilde{f}_n(\bm{\lambda})=1,
  \ \ \tilde{b}_n(\bm{\lambda})=-q^{-n}(1-cq^{n-1})(1-bd^{-1}q^{n-N}),
  \ \bm{\bar{\delta}}=(0,0,1,1),\!\\
  \text{(d)}&:
  \ \tilde{f}_n(\bm{\lambda})=-q^{-n}(1-cq^n)(1-bd^{-1}q^{n-N-1}),
  \ \ \tilde{b}_n(\bm{\lambda})=1,
  \ \bm{\bar{\delta}}=(0,0,-1,-1),\\
  \text{(e)}&:
  \ \tilde{f}_n(\bm{\lambda})=-q^{-n}(1-bq^n)(1-cd^{-1}q^{n-N-1}),
  \ \ \tilde{b}_n(\bm{\lambda})=1,
  \ \bm{\bar{\delta}}=(0,-1,0,-1),\\
  \text{(f)}&:
  \ \tilde{f}_n(\bm{\lambda})
  =\mathcal{E}_N(\bm{\lambda})-\mathcal{E}_n(\bm{\lambda}),
  \ \ \tilde{b}_n(\bm{\lambda})=1,
  \ \bm{\bar{\delta}}=(-1,0,0,-1).
\end{align}
For other cases, see Appendix \ref{app:rdQM}.

Then we can show that $B_1(x)$, $B_2(x)$, $D_1(x)$ and $D_2(x)$ satisfy
\begin{equation}
  B_1(x-1;\bm{\lambda})D_2(x;\bm{\lambda})
  +D_1(x+1;\bm{\lambda})B_2(x;\bm{\lambda})
  -\tilde{f}_0(\bm{\lambda})\tilde{b}_0(\bm{\lambda})
  =-B(x;\bm{\lambda})-D(x;\bm{\lambda}),
  \label{B1D2...}
\end{equation}
and the constants $\tilde{f}_n$ and $\tilde{b}_n$ satisfy
\begin{equation}
  \mathcal{E}_n(\bm{\lambda})
  =\tilde{f}_0(\bm{\lambda})\tilde{b}_0(\bm{\lambda})
  -\tilde{f}_n(\bm{\lambda})\tilde{b}_n(\bm{\lambda})
  \ \ (n\in\mathbb{Z}_{\geq 0}).
  \label{En=ft0bt0-..}
\end{equation}
The relation \eqref{B1D2...} gives another factorization of
$\widetilde{\mathcal{H}}(\bm{\lambda})$ \eqref{tHrdQM},
\begin{equation}
  \widetilde{\mathcal{H}}(\bm{\lambda})
  =-\tilde{\mathcal{B}}(\bm{\lambda})\tilde{\mathcal{F}}(\bm{\lambda})
  +\tilde{f}_0(\bm{\lambda})\tilde{b}_0(\bm{\lambda}).
  \label{tH=-tBtF..}
\end{equation}
Corresponding to this factorization \eqref{tH=-tBtF..},
we obtain the following relations.
\begin{mythm}
\label{thm:newfbsr_rdQM}
For the polynomials in \eqref{polyrdQMf}--\eqref{polyrdQMsi}
(except {\rm C} and {\rm $q$B}),
the following forward and backward shift relations hold for
$n\in\mathbb{Z}_{\geq 0}$,
\begin{align}
  \tilde{\mathcal{F}}(\bm{\lambda})\check{P}_n(x;\bm{\lambda})
  &=\tilde{f}_n(\bm{\lambda})\check{P}_n(x+s;\bm{\lambda}-\bm{\bar{\delta}}),
  \label{tFcPn=rdQM}\\
  \tilde{\mathcal{B}}(\bm{\lambda})
  \check{P}_n(x+s;\bm{\lambda}-\bm{\bar{\delta}})
  &=\tilde{b}_n(\bm{\lambda})\check{P}_n(x;\bm{\lambda}),
  \label{tBcPn=rdQM}
\end{align}
where $s$ is given by
\begin{equation}
  s=\left\{\begin{array}{ll}
  1&:\text{\rm H\,(a)(b),\,K\,(a),\,R\,(a)(b)(c),\,dH\,(a)(b)(c),
  \!dq$q$K\,(a)(b),\,$q$H\,(a)(b),}\\
  &\ \ \text{\rm $q$K\,(a),\,q$q$K\,(a)(b),\,a$q$K\,(a)(b),
  \!$q$R\,(a)(b)(c),\,d$q$H\,(a)(b)(c),}\\
  &\ \ \text{\rm d$q$K\,(a)(b),\,M\,(a),\,l$q$J\,(a),\,l$q$L\,(a),
  \!$q$M\,(a)(b),\,ASC$\II$\,(a),\,$q$C\,(a)}\\
  0&:\text{others}
  \end{array}\right..
\end{equation}
\end{mythm}
Proof:
It is sufficient to show \eqref{tFcPn=rdQM}, because \eqref{tHcPn=} and
\eqref{En=ft0bt0-..}--\eqref{tFcPn=rdQM} imply \eqref{tBcPn=rdQM}.
Taking $q$R (a) \eqref{B1B2D1D2:qR} as an example, let us prove
\eqref{tFcPn=rdQM}.
It is shown by direct calculation:
\begin{align*}
  &\quad\tilde{\mathcal{F}}(\bm{\lambda})\check{P}_n(x;\bm{\lambda})\\
  &=\frac{(1-a^{-1}dq^{x+1})(1-q^{x+1})}{(1-a^{-1}q)(1-dq^{2x+1})}
  {}_4\phi_3\Bigl(
  \genfrac{}{}{0pt}{}{q^{-n},\,abcd^{-1}q^{n-1},\,q^{-x},\,dq^x}
  {a,\,b,\,c}\Bigm|q\,;q\Bigr)\\
  &\quad-\frac{(1-aq^x)(1-dq^x)}{(aq^{-1}-1)(1-dq^{2x+1})}
  {}_4\phi_3\Bigl(
  \genfrac{}{}{0pt}{}{q^{-n},\,abcd^{-1}q^{n-1},\,q^{-x-1},\,dq^{x+1}}
  {a,\,b,\,c}\Bigm|q\,;q\Bigr)\\
  &=\frac{1}{(1-aq^{-1})(1-dq^{2x+1})}
  \sum_{k=0}^n\frac{(q^{-n},abcd^{-1}q^{n-1},q^{-x-1},dq^x\,;q)_k}
  {(a,b,c\,;q)_k}\frac{q^k}{(q\,;q)_k}\\
  &\qquad\times
  \bigl(-aq^{-1}(1-a^{-1}dq^{x+1})(-q^{x+1})(1-q^{-x+k-1})
  +(1-aq^x)(1-dq^{x+k})\bigr)\\
  &=\frac{1}{(1-aq^{-1})(1-dq^{2x+1})}
  \sum_{k=0}^n\frac{(q^{-n},abcd^{-1}q^{n-1},q^{-x-1},dq^x\,;q)_k}
  {(a,b,c\,;q)_k}\frac{q^k}{(q\,;q)_k}
  (1-aq^{k-1})(1-dq^{2x+1})\\
  &=\sum_{k=0}^n\frac{(q^{-n},abcd^{-1}q^{n-1},q^{-x-1},dq^x\,;q)_k}
  {(aq^{-1},b,c\,;q)_k}\frac{q^k}{(q\,;q)_k}\\
  &=\tilde{f}_n(\bm{\lambda})\check{P}_n(x+s;\bm{\lambda}-\bm{\bar{\delta}}).
\end{align*}
The other cases are proved in the same way.
\hfill\fbox{}

\medskip

\remark
Two formulas with $\bm{\bar{\delta}}$ and $-\bm{\bar{\delta}}$ are equivalent
by interchanging $\tilde{\mathcal{F}}$ and $\tilde{\mathcal{B}}$, e.g.
\eqref{tFcPn=rdQM} and \eqref{tBcPn=rdQM} for H\,(c) agree with
\eqref{tBcPn=rdQM} and \eqref{tFcPn=rdQM} for H\,(b) with the replacements
$a\to a+1$ and $b\to b-1$, respectively.
For C and $q$B, we do not have new factorization \eqref{tH=-tBtF..} and
another type of forward and backward shift relations
\eqref{tFcPn=rdQM}--\eqref{tBcPn=rdQM}.

\remark\label{xshift}\!\!\!
The relations \eqref{tFcPn=rdQM}--\eqref{tBcPn=rdQM} for twelve cases
((a) of H, K, R, dH, dq$q$K, $q$H, $q$K, q$q$K, a$q$K, $q$R, d$q$H, d$q$K,
which have $\tilde{f}_n=1$, $\tilde{b}_n(\bm{\lambda})=
\mathcal{E}_{N+1}(\bm{\lambda})-\mathcal{E}_n(\bm{\lambda})$, $s=1$ and
$D_1(0;\bm{\lambda})=B_1(N;\bm{\lambda})=0$) were given in \cite{os40} and
they were called forward and backward $x$-shift relations.
By considering $e^{-\partial}\tilde{\mathcal{F}}(\bm{\lambda})$ and
$\tilde{\mathcal{B}}(\bm{\lambda})e^{\partial}$, the above results
\eqref{tH=-tBtF..} and \eqref{tFcPn=rdQM}--\eqref{tBcPn=rdQM} with $s=1$ are
rewritten as
\begin{align}
  \widetilde{\mathcal{H}}(\bm{\lambda})
  =-\bigl(\tilde{\mathcal{B}}(\bm{\lambda})e^{\partial}\bigr)
  &\bigl(e^{-\partial}\tilde{\mathcal{F}}(\bm{\lambda})\bigr)
  +\tilde{f}_0(\bm{\lambda})\tilde{b}_0(\bm{\lambda}),\\
  \bigl(e^{-\partial}\tilde{\mathcal{F}}(\bm{\lambda})\bigr)
  \check{P}_n(x;\bm{\lambda})
  &=\tilde{f}_n(\bm{\lambda})\check{P}_n(x;\bm{\lambda}-\bm{\bar{\delta}}),\\
  \bigl(\tilde{\mathcal{B}}(\bm{\lambda})e^{\partial}\bigr)
  \check{P}_n(x;\bm{\lambda}-\bm{\bar{\delta}})
  &=\tilde{b}_n(\bm{\lambda})\check{P}_n(x;\bm{\lambda}).
\end{align}
That is, $x$ is not shifted. As an identity of polynomial, the $x$-shift is not
essential. However, we think that this $x$-shift has important implications in
the state-adding Darboux transformation for the finite rdQM systems
\cite{os40,addrdQM}.

\remark\label{AWqR}\!\!\!
%AW and $q$R polynomials are related as \cite{kls,casoidrdqm}
AW and $q$R polynomials are related as \cite{kls}
\begin{align}
  &e^{ix^{\text{AW}}}=d^{\frac12}q^{x^{\text{$q$R}}},\quad
  (a_1,a_2,a_3,a_4)=(ad^{-\frac12},bd^{-\frac12},cd^{-\frac12},d^{\frac12}),\n
  &\check{P}^{\text{AW}}_n(x^{\text{AW}};\bm{\lambda}^{\text{AW}})
  =d^{-\frac{n}{2}}(a,b,c\,;q)_n
  \check{P}^{\text{$q$R}}_n(x^{\text{$q$R}};\bm{\lambda}^{\text{$q$R}}).
\end{align}
For the $(j,k)=(1,4)$ case in \eqref{V1V2:AW}, the operators
$\tilde{\mathcal{F}}$ and $\tilde{\mathcal{B}}$ for AW are related to those
for $q$R (a) \eqref{B1B2D1D2:qR} as
\begin{align}
  e^{\frac{\gamma}{2}p}\tilde{\mathcal{F}}^{\text{AW}}(\bm{\lambda}^{\text{AW}})
  &=-(q^{-N-1}-1)\tilde{\mathcal{F}}^{\text{$q$R}}
  (\bm{\lambda}^{\text{$q$R}}),\n
  \tilde{\mathcal{B}}^{\text{AW}}(\bm{\lambda}^{\text{AW}})
  e^{-\frac{\gamma}{2}p}
  &=(q^{-N-1}-1)^{-1}
  \tilde{\mathcal{B}}^{\text{$q$R}}(\bm{\lambda}^{\text{$q$R}}).
\end{align}
These extra factors $e^{\pm\frac{\gamma}{2}p}$ give the property in
Remark\,\ref{*inv}.
Similarly AW with $(j,k)=(2,4)$, $(3,4)$, $(1,2)$, $(1,3)$ and $(2,4)$ cases
correspond to $q$R (b), (c), (d), (e) and (f), respectively.

\remark\label{simodoki:rdQM}\!\!\!
We can show that
\begin{align}
  &B_1(x-s;\bm{\lambda})B_2(x-s+1;\bm{\lambda})
  =B(x;\bm{\lambda}-\bm{\bar{\delta}}),
  \ \ D_1(x-s+1;\bm{\lambda})D_2(x-s;\bm{\lambda})
  =D(x;\bm{\lambda}-\bm{\bar{\delta}}),\n
  &B_1(x-s;\bm{\lambda})D_2(x-s+1;\bm{\lambda})
  +D_1(x-s+1;\bm{\lambda})B_2(x-s;\bm{\lambda})
  -\tilde{f}_0(\bm{\lambda})\tilde{b}_0(\bm{\lambda})\\
  &=B_1(x-1;\bm{\lambda}-\bm{\bar{\delta}})
  D_2(x;\bm{\lambda}-\bm{\bar{\delta}})
  +D_1(x+1;\bm{\lambda}-\bm{\bar{\delta}})
  B_2(x;\bm{\lambda}-\bm{\bar{\delta}})
  -\tilde{f}_0(\bm{\lambda}-\bm{\bar{\delta}})
  \tilde{b}_0(\bm{\lambda}-\bm{\bar{\delta}}),
 \nonumber
\end{align}
which imply
\begin{equation}
  \tilde{\mathcal{F}}(\bm{\lambda})\tilde{\mathcal{B}}(\bm{\lambda})
  \Bigm|_{x\to x-s}
  -\tilde{f}_0(\bm{\lambda})\tilde{b}_0(\bm{\lambda})
  =\tilde{\mathcal{B}}(\bm{\lambda}-\bm{\bar{\delta}})
  \tilde{\mathcal{F}}(\bm{\lambda}-\bm{\bar{\delta}})
  -\tilde{f}_0(\bm{\lambda}-\bm{\bar{\delta}})
  \tilde{b}_0(\bm{\lambda}-\bm{\bar{\delta}}).
  \label{tFtB-..:rdQM}
\end{equation}

%%%%%%%%%%%%%%%%%%%%%%%%%%%%%%%%%%%%%%%%%%%%%%%%%%%%%
%                                                   %
% 4.4 Polynomials in rdQMJ systems                  %
%                                                   %
%%%%%%%%%%%%%%%%%%%%%%%%%%%%%%%%%%%%%%%%%%%%%%%%%%%%%
\subsection{Polynomials in rdQMJ systems}
\label{sec:newrdQMJ}
\addtocounter{mybangouII}{1}
\setcounter{myrembangou}{0}

For the rdQMJ systems described by the polynomials \eqref{polyrdQMJ}
(except d$q$He$\I$, d$q$He$\II$ and SW),
let us define the operators $\tilde{\mathcal{F}}^{\text{J}}(\bm{\lambda})$ and
$\tilde{\mathcal{B}}^{\text{J}}(\bm{\lambda})$ as follows:
\begin{align}
  \tilde{\mathcal{F}}^{\text{J}}(\bm{\lambda})&\eqdef
  D^{\text{J}}_1(q\eta;\bm{\lambda})+B^{\text{J}}_1(\eta;\bm{\lambda})
  q^{\eta\frac{d}{d\eta}},
  \label{tFrdQMJ}\\
  \tilde{\mathcal{B}}^{\text{J}}(\bm{\lambda})&\eqdef
  B^{\text{J}}_2(\eta;\bm{\lambda})+D^{\text{J}}_2(\eta;\bm{\lambda})
  q^{-\eta\frac{d}{d\eta}},
  \label{tBrdQMJ}
\end{align}
where the potential functions $B^{\text{J}}_1(\eta;\bm{\lambda}),
B^{\text{J}}_2(\eta;\bm{\lambda}),D^{\text{J}}_1(\eta;\bm{\lambda})$ and
$D^{\text{J}}_2(\eta;\bm{\lambda})$ satisfy
\begin{equation}
  B^{\text{J}}(\eta;\bm{\lambda})=
  B^{\text{J}}_1(\eta;\bm{\lambda})B^{\text{J}}_2(\eta;\bm{\lambda}),\quad
  D^{\text{J}}(\eta;\bm{\lambda})=
  D^{\text{J}}_1(\eta;\bm{\lambda})D^{\text{J}}_2(\eta;\bm{\lambda}).
  \label{BJ=BJ1BJ2}
\end{equation}
For b$q$J case, their explicit forms are given by
\begin{align}
%%%%%%%%%%%%%%%%%%%%%%%%%%%%%%%%
% big $q$-Jacobi (b$q$J)       %
%%%%%%%%%%%%%%%%%%%%%%%%%%%%%%%%
  \text{(a)}&:
  \ B^{\text{J}}_1(\eta;\bm{\lambda})=\frac{\eta^{-1}a(1-\eta)}{1-a},
  \ \ B^{\text{J}}_2(\eta;\bm{\lambda})=(1-a)\eta^{-1}q(b\eta-c),\n
  &\ \ \ D^{\text{J}}_1(\eta;\bm{\lambda})=\frac{\eta^{-1}(aq-\eta)}{a-1},
  \ \ D^{\text{J}}_2(\eta;\bm{\lambda})=(a-1)\eta^{-1}(\eta-cq),
  \label{B1JB2JD1JD2J:bqJ}\\
  \text{(b)}&:
  \ B^{\text{J}}_1(\eta;\bm{\lambda})=\frac{\eta^{-1}(1-\eta)}{c^{-1}-1},
  \ \ B^{\text{J}}_2(\eta;\bm{\lambda})=(c^{-1}-1)\eta^{-1}aq(b\eta-c),\n
  &\ \ \ D^{\text{J}}_1(\eta;\bm{\lambda})=\frac{\eta^{-1}(\eta-cq)}{1-c},
  \ \ D^{\text{J}}_2(\eta;\bm{\lambda})=(1-c)\eta^{-1}(aq-\eta),\\
  \text{(c)}&:
  \ B^{\text{J}}_1(\eta;\bm{\lambda})=(c^{-1}q^{-1}-1)\eta^{-1}aq(b\eta-c),
  \ \ B^{\text{J}}_2(\eta;\bm{\lambda})=\frac{\eta^{-1}(1-\eta)}
  {c^{-1}q^{-1}-1},\n
  &\ \ \ D^{\text{J}}_1(\eta;\bm{\lambda})=(1-cq)\eta^{-1}(aq-\eta),
  \ \ D^{\text{J}}_2(\eta;\bm{\lambda})=\frac{\eta^{-1}(\eta-cq)}{1-cq},\\
  \text{(d)}&:
  \ B^{\text{J}}_1(\eta;\bm{\lambda})=(1-aq)\eta^{-1}(b\eta-c),
  \ \ B^{\text{J}}_2(\eta;\bm{\lambda})=\frac{\eta^{-1}aq(1-\eta)}{1-aq},\n
  &\ \ \ D^{\text{J}}_1(\eta;\bm{\lambda})=(aq-1)\eta^{-1}(\eta-cq),
  \ \ D^{\text{J}}_2(\eta;\bm{\lambda})=\frac{\eta^{-1}(aq-\eta)}{aq-1},
\end{align}
and the constants $\tilde{f}^{\text{J}}_n(\bm{\lambda})$,
$\tilde{b}^{\text{J}}_n(\bm{\lambda})$ and $\bm{\bar{\delta}}$ are given by
\begin{align}
%%%%%%%%%%%%%%%%%%%%%%%%%%%%%%%%
% big $q$-Jacobi (b$q$J)       %
%%%%%%%%%%%%%%%%%%%%%%%%%%%%%%%%
  \text{(a)}&:
  \ \tilde{f}^{\text{J}}_n(\bm{\lambda})=1,
  \ \ \tilde{b}^{\text{J}}_n(\bm{\lambda})=-q^{-n}(1-aq^n)(1-bq^{n+1}),
  \ \ \bm{\bar{\delta}}=(1,-1,0),\\
  \text{(b)}&:
  \ \tilde{f}^{\text{J}}_n(\bm{\lambda})=1,
  \ \ \tilde{b}^{\text{J}}_n(\bm{\lambda})=-q^{-n}(1-cq^n)(1-abc^{-1}q^{n+1}),
  \ \ \bm{\bar{\delta}}=(0,0,1),\\
  \text{(c)}&:
  \ \tilde{f}^{\text{J}}_n(\bm{\lambda})=-q^{-n}(1-cq^{n+1})(1-abc^{-1}q^n),
  \ \ \tilde{b}^{\text{J}}_n(\bm{\lambda})=1,
  \ \ \bm{\bar{\delta}}=(0,0,-1),\\
  \text{(d)}&:
  \ \tilde{f}^{\text{J}}_n(\bm{\lambda})=-q^{-n}(1-aq^{n+1})(1-bq^n),
  \ \ \tilde{b}^{\text{J}}_n(\bm{\lambda})=1,
  \ \ \bm{\bar{\delta}}=(-1,1,0).
\end{align}
For other cases, see Appendix \ref{app:rdQMJ}.

Then we can show that $B^{\text{J}}_1(\eta)$, $B^{\text{J}}_2(\eta)$,
$D^{\text{J}}_1(\eta)$ and $D^{\text{J}}_2(\eta)$ satisfy
\begin{equation}
  B^{\text{J}}_1(q^{-1}\eta;\bm{\lambda})D^{\text{J}}_2(\eta;\bm{\lambda})
  +D^{\text{J}}_1(q\eta;\bm{\lambda})B^{\text{J}}_2(\eta;\bm{\lambda})
  -\tilde{f}^{\text{J}}_0(\bm{\lambda})\tilde{b}^{\text{J}}_0(\bm{\lambda})
  =-B^{\text{J}}(\eta;\bm{\lambda})-D^{\text{J}}(\eta;\bm{\lambda}),
  \label{BJ1BJ2=...}
\end{equation}
and the constants $\tilde{f}^{\text{J}}_n$ and $\tilde{b}^{\text{J}}_n$ satisfy
\begin{equation}
  \mathcal{E}_n(\bm{\lambda})
  =\tilde{f}^{\text{J}}_0(\bm{\lambda})\tilde{b}^{\text{J}}_0(\bm{\lambda})
  -\tilde{f}^{\text{J}}_n(\bm{\lambda})\tilde{b}^{\text{J}}_n(\bm{\lambda})
  \ \ (n\in\mathbb{Z}_{\geq 0}).
  \label{En=ftJ0btJ0-..}
\end{equation}
The relations \eqref{BJ1BJ2=...} give other factorizations of
$\widetilde{\mathcal{H}}(\bm{\lambda})$ \eqref{tHrdQMJ},
\begin{equation}
  \widetilde{\mathcal{H}}^{\text{J}}(\bm{\lambda})
  =-\tilde{\mathcal{B}}^{\text{J}}(\bm{\lambda})
  \tilde{\mathcal{F}}^{\text{J}}(\bm{\lambda})
  +\tilde{f}^{\text{J}}_0(\bm{\lambda})\tilde{b}^{\text{J}}_0(\bm{\lambda}).
  \label{tHJ=-tBJtFJ..}
\end{equation}
Corresponding to this factorization \eqref{tHJ=-tBJtFJ..},
we obtain the following relations.
\begin{mythm}
\label{thm:newfbsr_rdQMJ}
For the polynomials in \eqref{polyrdQMJ} (except {\rm d$q$He$\rm\I$,
d$q$He$\rm\II$} and {\rm SW}), the following forward and backward shift
relations hold for $n\in\mathbb{Z}_{\geq 0}$,
\begin{align}
  \tilde{\mathcal{F}}^{\text{\rm J}}(\bm{\lambda})P_n(\eta;\bm{\lambda})
  &=\tilde{f}^{\text{\rm J}}_n(\bm{\lambda})
  P_n(r'\eta;\bm{\lambda}-\bm{\bar{\delta}}),
  \label{tFJPn=rdQMJ}\\
  \tilde{\mathcal{B}}^{\text{\rm J}}(\bm{\lambda})
  P_n(r'\eta;\bm{\lambda}-\bm{\bar{\delta}})
  &=\tilde{b}^{\text{\rm J}}_n(\bm{\lambda})P_n(\eta;\bm{\lambda}),
  \label{tBJPn=rdQMJ}
\end{align}
where $r'$ is given by
\begin{equation}
  r'=\left\{
  \begin{array}{ll}
  q&:\text{\rm b$q$J\,(c)(d),\,b$q$L\,(c)(d),\,ASC$\I$\,(a),\,$q$L\,(b)}
  \\[2pt]
  1&:\text{\rm b$q$J\,(a)(b),\,b$q$L\,(a)(b),\,ASC$\I$\,(b),\,$q$L\,(a)}
  \end{array}\right..
\end{equation}
\end{mythm}
Proof:
It is sufficient to show \eqref{tFJPn=rdQMJ}, because \eqref{tHJPn=} and
\eqref{En=ftJ0btJ0-..}--\eqref{tFJPn=rdQMJ} imply \eqref{tBJPn=rdQMJ}.
Taking b$q$J (a) \eqref{B1JB2JD1JD2J:bqJ} as an example, let us prove
\eqref{tFJPn=rdQMJ}.
It is shown by direct calculation:
\begin{align*}
  &\quad\tilde{\mathcal{F}}^{\text{J}}(\bm{\lambda})P_n(\eta;\bm{\lambda})\\
  &=\frac{\eta^{-1}(a-\eta)}{a-1}
  {}_3\phi_2\Bigl(\genfrac{}{}{0pt}{}{q^{-n},\,abq^{n+1},\,\eta}
  {aq,\,cq}\Bigm|q\,;q\Bigr)
  +\frac{\eta^{-1}a(1-\eta)}{1-a}
  {}_3\phi_2\Bigl(\genfrac{}{}{0pt}{}{q^{-n},\,abq^{n+1},\,q\eta}
  {aq,\,cq}\Bigm|q\,;q\Bigr)\\
  &=\frac{\eta^{-1}}{1-a}
  \sum_{k=0}^n\frac{(q^{-n},abq^{n+1},\eta\,;q)_k}
  {(aq,cq\,;q)_k}\frac{q^k}{(q\,;q)_k}
  \bigl(-(a-\eta)+a(1-\eta q^k)\bigr)\\
  &=\frac{\eta^{-1}}{1-a}
  \sum_{k=0}^n\frac{(q^{-n},abq^{n+1},\eta\,;q)_k}
  {(aq,cq\,;q)_k}\frac{q^k}{(q\,;q)_k}
  (1-aq^k)\eta\\
  &=\sum_{k=0}^n\frac{(q^{-n},abq^{n+1},\eta\,;q)_k}
  {(a,cq\,;q)_k}\frac{q^k}{(q\,;q)_k}\\
  &=\tilde{f}^{\text{J}}_n(\bm{\lambda})
  P_n(r'\eta;\bm{\lambda}-\bm{\bar{\delta}}).
\end{align*}
The other cases are proved in the same way.
\hfill\fbox{}

\medskip

\remark
Two formulas with $\bm{\bar{\delta}}$ and $-\bm{\bar{\delta}}$ are equivalent
by interchanging $\tilde{\mathcal{F}}^{\text{J}}$ and
$\tilde{\mathcal{B}}^{\text{J}}$, e.g.
\eqref{tFJPn=rdQMJ} and \eqref{tBJPn=rdQMJ} for b$q$J\,(d) agree with
\eqref{tBJPn=rdQMJ} and \eqref{tFJPn=rdQMJ} for b$q$J\,(a) with the replacements
$a\to aq$ and $b\to bq^{-1}$, respectively.
For d$q$He$\I$, d$q$He$\II$ and SW, we do not have new factorization
\eqref{tHJ=-tBJtFJ..} and another type of forward and backward shift relations
\eqref{tFJPn=rdQMJ}--\eqref{tBJPn=rdQMJ}.

\remark\label{xshiftJ}\!\!\!
As in Remark\,\ref{xshift}, by considering
$q^{-\eta\frac{d}{d\eta}}\tilde{\mathcal{F}}^{\text{J}}(\bm{\lambda})$ and
$\tilde{\mathcal{B}}^{\text{J}}(\bm{\lambda})q^{\eta\frac{d}{d\eta}}$,
the above results \eqref{tHJ=-tBJtFJ..} and
\eqref{tFJPn=rdQMJ}--\eqref{tBJPn=rdQMJ} with $r'=q$ are rewritten as
\begin{align}
  \widetilde{\mathcal{H}}^{\text{J}}(\bm{\lambda})
  =-\bigl(\tilde{\mathcal{B}}^{\text{J}}(\bm{\lambda})
  q^{\eta\frac{d}{d\eta}}\bigr)
  &\bigl(q^{-\eta\frac{d}{d\eta}}
  \tilde{\mathcal{F}}^{\text{J}}(\bm{\lambda})\bigr)
  +\tilde{f}^{\text{J}}_0(\bm{\lambda})\tilde{b}^{\text{J}}_0(\bm{\lambda}),\\
  \bigl(q^{-\eta\frac{d}{d\eta}}
  \tilde{\mathcal{F}}^{\text{J}}(\bm{\lambda})\bigr)P_n(\eta;\bm{\lambda})
  &=\tilde{f}^{\text{J}}_n(\bm{\lambda})
  P_n(\eta;\bm{\lambda}-\bm{\bar{\delta}}),\\
  \bigl(\tilde{\mathcal{B}}^{\text{J}}(\bm{\lambda})
  q^{\eta\frac{d}{d\eta}}\bigr)P_n(\eta;\bm{\lambda}-\bm{\bar{\delta}})
  &=\tilde{b}^{\text{J}}_n(\bm{\lambda})P_n(\eta;\bm{\lambda}).
\end{align}
That is, $\eta$ is not $q$-shifted. As an identity of polynomial,
the $q$-shift of $\eta$ is not essential.

\remark\label{simodoki:rdQMJ}\!\!\!
We can show that
\begin{align}
  &B^{\text{J}}_1(r^{\prime\,-1}\eta;\bm{\lambda})
  B^{\text{J}}_2(qr^{\prime\,-1}\eta;\bm{\lambda})
  =B^{\text{J}}(\eta;\bm{\lambda}-\bm{\bar{\delta}}),
  \ \ D^{\text{J}}_1(qr^{\prime\,-1}\eta;\bm{\lambda})
  D^{\text{J}}_2(r^{\prime\,-1}\eta;\bm{\lambda})
  =D^{\text{J}}(\eta;\bm{\lambda}-\bm{\bar{\delta}}),\n
  &B^{\text{J}}_1(r^{\prime\,-1}\eta;\bm{\lambda})
  D^{\text{J}}_2(qr^{\prime\,-1}\eta;\bm{\lambda})
  +D^{\text{J}}_1(qr^{\prime\,-1}\eta;\bm{\lambda})
  B^{\text{J}}_2(r^{\prime\,-1}\eta;\bm{\lambda})
  -\tilde{f}^{\text{J}}_0(\bm{\lambda})\tilde{b}^{\text{J}}_0(\bm{\lambda})\\
  &=B^{\text{J}}_1(q^{-1}\eta;\bm{\lambda}-\bm{\bar{\delta}})
  D^{\text{J}}_2(\eta;\bm{\lambda}-\bm{\bar{\delta}})
  +D^{\text{J}}_1(q\eta;\bm{\lambda}-\bm{\bar{\delta}})
  B^{\text{J}}_2(\eta;\bm{\lambda}-\bm{\bar{\delta}})
  -\tilde{f}^{\text{J}}_0(\bm{\lambda}-\bm{\bar{\delta}})
  \tilde{b}^{\text{J}}_0(\bm{\lambda}-\bm{\bar{\delta}}),
 \nonumber
\end{align}
which imply
\begin{equation}
  \tilde{\mathcal{F}}^{\text{J}}(\bm{\lambda})
  \tilde{\mathcal{B}}^{\text{J}}(\bm{\lambda})
  \Bigm|_{\eta\to r^{\prime\,-1}\eta}
  -\tilde{f}^{\text{J}}_0(\bm{\lambda})\tilde{b}^{\text{J}}_0(\bm{\lambda})
  =\tilde{\mathcal{B}}^{\text{J}}(\bm{\lambda}-\bm{\bar{\delta}})
  \tilde{\mathcal{F}}^{\text{J}}(\bm{\lambda}-\bm{\bar{\delta}})
  -\tilde{f}^{\text{J}}_0(\bm{\lambda}-\bm{\bar{\delta}})
  \tilde{b}^{\text{J}}_0(\bm{\lambda}-\bm{\bar{\delta}}).
  \label{tFtB-..:rdQMJ}
\end{equation}

%%%%%%%%%%%%%%%%%%%%%%%%%%%%%%%%%%%%%%%%%%%%%%%%%%%%%%%%%%%%%%%
%                                                             %
%  5. Summary and Comments                                    %
%                                                             %
%%%%%%%%%%%%%%%%%%%%%%%%%%%%%%%%%%%%%%%%%%%%%%%%%%%%%%%%%%%%%%%
\section{Summary and Comments}
\label{sec:summary}

The orthogonal polynomials in the Askey scheme satisfy second order
differential or difference equations (Theorem\,\ref{thm:sabuneq}) and
we study them by using quantum mechanical formulation (oQM, idQM, rdQM, rdQMJ).
The forward and backward shift relations are their basic properties
(Theorem\,\ref{thm:fbsr}, \ref{thm:fbsrJ}), in which the degree $n$ and the
parameters $\bm{\lambda}$ are shifted.
They are based on the factorizations of the differential or difference
operators $\widetilde{\mathcal{H}}$ \eqref{tH=BF} and
$\widetilde{\mathcal{H}}^{\text{J}}$ \eqref{tHJ=BJFJ}.
Although Theorem\,\ref{thm:sabuneq}, \ref{thm:fbsr} and \ref{thm:fbsrJ} are
well known results, their formulas in quantum mechanical formulation are
expressed neatly and systematically in a universal form.
Motivated by the recently found forward and backward $x$-shift relations
\cite{os40}, in which the coordinate $x$ and parameters $\bm{\lambda}$ are
shifted, we have tried to find new forward and backward relations.
We have found new factorizations of $\widetilde{\mathcal{H}}$
\eqref{tH=tBtF..oQM}, \eqref{tH=tBtF..idQM}, \eqref{tH=-tBtF..}
and $\widetilde{\mathcal{H}}^{\text{J}}$ \eqref{tHJ=-tBJtFJ..}, and based on
them, we have obtained another type of forward and backward shift relations
(Theorem\,\ref{thm:newfbsr_oQM}, \ref{thm:newfbsr_idQM},
\ref{thm:newfbsr_rdQM}, \ref{thm:newfbsr_rdQMJ}).
While some of these results may be found in the literature, this is the first
comprehensive study.
In these forward and backward shift relations except for some cases of
rdQM and rdQMJ, only the parameters $\bm{\lambda}$ are shifted.
As an identity of polynomial, the $x$-shift (or $q$-shift of $\eta$) is not
essential (Remark\,\ref{xshift}, \ref{xshiftJ}).

The forward and backward shift relations are related to the shape invariance
property of quantum mechanical systems \cite{os12,os13,os24,os34}.
It is an interesting problem to investigate the quantum mechanical implications
of the another type of forward and backward shift relations obtained in this
paper
(cf.\ Remark\,\ref{simodoki:idQM}, \ref{simodoki:rdQM}, \ref{simodoki:rdQMJ}).
Especially the twelve finite rdQM cases in Remark\,\ref{xshift} are interesting.
In these cases, we think that the $x$-shift has important implications related
to the state-adding Darboux transformations.
The state-adding Darboux transformations for finite rdQM systems were studied
in \cite{os40,addrdQM}. For one-step transformation, the range of $x$,
$\{0,1,\ldots,N\}$, is extended to $\{-1,0,1,\ldots,N\}$, and the parameter $N$
($=(\text{size of the Hamiltonian})-1$)
is shifted to $N+1$, and the deformed potential functions contain the factors
$B(x+1;\bm{\lambda}-\bm{\bar{\delta}})$ and
$D(x+1;\bm{\lambda}-\bm{\bar{\delta}})$, where $\bm{\bar{\delta}}$ is certain
shift of parameters (the component of $\bm{\bar{\delta}}$ corresponding to $N$
is $-1$).
The boundary conditions $B(x)=0$ for $x=x_{\text{max}}=N$ and $D(x)=0$ for
$x=x_{\text{min}}=0$ are inherited by the deformed system, because we have
$B(x+1;\bm{\lambda}-\bm{\bar{\delta}})=0$ for $x=x_{\text{max}}=N$ and
$D(x+1;\bm{\lambda}-\bm{\bar{\delta}})=0$ for $x=x_{\text{min}}=-1$.
Thus the $x$-shift, $x\to x+1$, is important.
The $M$ added eigenvectors of the deformed systems by $M$-step state-adding
Darboux transformations are obtained explicitly in \cite{addrdQM}.
But they were found through a very technical trial and error process, and
a better derivation is desired.
For the state-adding Darboux transformations in oQM, the added eigenfunctions
are expressed neatly in terms of the Wronskian and seed solutions.
We believe that this is also the case for finite rdQM, namely, the added
eigenvectors are expressed neatly in terms of the Casoratian and seed solutions,
and that the $x$-shift relations in Remark\,\ref{xshift} play an important role.
We hope this topic will be successfully solved and we will be able to report
the results somewhere.

%The case-(1) multi-indexed orthogonal polynomials are constructed for
%R and $q$R \cite{os26}, W and AW \cite{os27}, M, l$q$J and l$q$L \cite{os35},
%cH and MP \cite{idQMcH}, and they have shape invariant property, namely,
The case-(1) multi-indexed orthogonal polynomials are constructed for
R, $q$R, W, AW, M, l$q$J, l$q$L, cH and MP,
and they have shape invariant property, namely,
satisfy the forward and backward shift relations like Theorem\,\ref{thm:fbsr}.
It is an interesting problem to investigate whether these multi-indexed
polynomials satisfy another type of forward and backward shift relations such as
Theorem\,\ref{thm:newfbsr_idQM} and \ref{thm:newfbsr_rdQM}.

%%%%%%%%%%%%%%%%%%%%%%%%%%%%%%%%%%%%%%%%%%%%%%%%%%%%%%%%%%%%%%%
%                                                             %
%  Acknowledgments                                            %
%                                                             %
%%%%%%%%%%%%%%%%%%%%%%%%%%%%%%%%%%%%%%%%%%%%%%%%%%%%%%%%%%%%%%%
\section*{Acknowledgements}

This work is supported by JSPS KAKENHI Grant Number JP19K03667.
%I thank Satoshi Tsujimoto for informing me of the paper \cite{km89}.

\bigskip
\appendix
%\renewcommand{\theequation}{\Alph{section}.\arabic{equation}}
%%%%%%%%%%%%%%%%%%%%%%%%%%%%%%%%%%%%%%%%%%%%%%%%%%%%%%%%%%%%%%%
%                                                             %
%  A. Data for \S\,\ref{sec:newfbsr}                          %
%                                                             %
%%%%%%%%%%%%%%%%%%%%%%%%%%%%%%%%%%%%%%%%%%%%%%%%%%%%%%%%%%%%%%%
\section{Data for \S\,\ref{sec:newfbsr}}
\label{app:data}

We give the data for the another type of forward and backward shift relations in
\S\,\ref{sec:newfbsr}.

%%%%%%%%%%%%%%%%%%%%%%%%%%%%%%%%%%%%%%%%%%%%%%%%%%%%%
%                                                   %
% A.1 Data for \S\,\ref{sec:newoQM}}                %
%                                                   %
%%%%%%%%%%%%%%%%%%%%%%%%%%%%%%%%%%%%%%%%%%%%%%%%%%%%%
\subsection{Data for \S\,\ref{sec:newoQM}}
\label{app:oQM}

We present explicit forms of $\tilde{\mathcal{F}}$, $\tilde{\mathcal{B}}$,
$\tilde{f}_n$, $\tilde{b}_n$ and $\bm{\bar{\delta}}$ in \S\,\ref{sec:newoQM}.

\noindent
Operators $\tilde{\mathcal{F}}(\bm{\lambda})$ and
$\tilde{\mathcal{B}}(\bm{\lambda})$:
\begin{align}
%%%%%%%%%%%%%%%%%%%%%%%%%%%%%%%%
% Laguerre (L)                 %
%%%%%%%%%%%%%%%%%%%%%%%%%%%%%%%%
  \text{L}&:\text{(a)}:
  \ \tilde{\mathcal{F}}(\bm{\lambda})\eqdef\frac12x\frac{d}{dx}+g-\frac12
  \ \Bigl(=\eta\frac{d}{d\eta}+g-\frac12\Bigr),
  \ \ \tilde{\mathcal{B}}(\bm{\lambda})\eqdef
  -\frac12\frac{1}{x}\frac{d}{dx}+1\ \Bigl(=-\frac{d}{d\eta}+1\Bigr),
  \label{tBoQM:L(a)}\\
  &\ \ \,\text{(b)}:
  \ \tilde{\mathcal{F}}(\bm{\lambda})\eqdef-\frac12\frac{1}{x}\frac{d}{dx}+1
  \ \Bigl(=-\frac{d}{d\eta}+1\Bigr),
  \ \ \tilde{\mathcal{B}}(\bm{\lambda})\eqdef
  \frac12x\frac{d}{dx}+g+\frac12\ \Bigl(=\eta\frac{d}{d\eta}+g+\frac12\Bigr),
  \label{tBoQM:L(b)}\\
%
%%%%%%%%%%%%%%%%%%%%%%%%%%%%%%%%
% pseudo Jacobi (pJ)           %
%%%%%%%%%%%%%%%%%%%%%%%%%%%%%%%%
  \text{pJ}&:\text{(a)}:
  \ \tilde{\mathcal{F}}(\bm{\lambda})\eqdef
  \Bigl(\tanh x+\frac{i}{\cosh x}\Bigr)\frac{d}{dx}-h-\frac12-i\mu
  \ \Bigl(=(\eta+i)\frac{d}{d\eta}-h-\frac12-i\mu\Bigr),\n
  &\phantom{:\text{(a)}:}\ \ \ \tilde{\mathcal{B}}(\bm{\lambda})\eqdef
  \Bigl(\tanh x-\frac{i}{\cosh x}\Bigr)\frac{d}{dx}-h+\frac12+i\mu
  \ \Bigl(=(\eta-i)\frac{d}{d\eta}-h+\frac12+i\mu\Bigr),
  \label{tBoQM:pJ(a)}\\
  &\ \ \,\text{(b)}:
  \ \tilde{\mathcal{F}}(\bm{\lambda})\eqdef
  \Bigl(\tanh x-\frac{i}{\cosh x}\Bigr)\frac{d}{dx}-h-\frac12+i\mu
  \ \Bigl(=(\eta-i)\frac{d}{d\eta}-h-\frac12+i\mu\Bigr),\n
  &\phantom{:\text{(b)}:}\ \ \ \tilde{\mathcal{B}}(\bm{\lambda})\eqdef
  \Bigl(\tanh x+\frac{i}{\cosh x}\Bigr)\frac{d}{dx}-h+\frac12-i\mu
  \ \Bigl(=(\eta+i)\frac{d}{d\eta}-h+\frac12-i\mu\Bigr).
  \label{tBoQM:pJ(b)}
\end{align}
Constants $\tilde{f}_n(\bm{\lambda})$, $\tilde{b}_n(\bm{\lambda})$ and
$\bm{\bar{\delta}}$:
\begin{align}
%%%%%%%%%%%%%%%%%%%%%%%%%%%%%%%%
% Laguerre (L)                 %
%%%%%%%%%%%%%%%%%%%%%%%%%%%%%%%%
  \text{L}&:\text{(a)}:
  \ \tilde{f}_n(\bm{\lambda})=n+g-\tfrac12,
  \ \ \tilde{b}_n(\bm{\lambda})=1,
  \ \ \bm{\bar{\delta}}=1,\\
  &\ \ \,\text{(b)}:
  \ \tilde{f}_n(\bm{\lambda})=1,
  \ \ \tilde{b}_n(\bm{\lambda})=n+g+\tfrac12,
  \ \ \bm{\bar{\delta}}=-1,\\
%
%%%%%%%%%%%%%%%%%%%%%%%%%%%%%%%%
% pseudo Jacobi (pJ)           %
%%%%%%%%%%%%%%%%%%%%%%%%%%%%%%%%
  \text{pJ}&:\text{(a)}:
  \ \tilde{f}_n(\bm{\lambda})=n-h-\tfrac12-i\mu,
  \ \ \tilde{b}_n(\bm{\lambda})=n-h+\tfrac12+i\mu,
  \ \ \bm{\bar{\delta}}=(0,i),\\
  &\ \ \,\text{(b)}:
  \ \tilde{f}_n(\bm{\lambda})=n-h-\tfrac12+i\mu,
  \ \ \tilde{b}_n(\bm{\lambda})=n-h+\tfrac12-i\mu,
  \ \ \bm{\bar{\delta}}=(0.-i).
\end{align}
We remark that the second components of $\bm{\bar{\delta}}$ for pJ are
unphysical values.

%%%%%%%%%%%%%%%%%%%%%%%%%%%%%%%%%%%%%%%%%%%%%%%%%%%%%
%                                                   %
% A.2 Data for \S\,\ref{sec:newidQM}                %
%                                                   %
%%%%%%%%%%%%%%%%%%%%%%%%%%%%%%%%%%%%%%%%%%%%%%%%%%%%%
\subsection{Data for \S\,\ref{sec:newidQM}}
\label{app:idQM}

We present explicit forms of $V_1(x)$, $\tilde{f}_n$, $\tilde{b}_n$ and
$\bm{\bar{\delta}}$ in \S\,\ref{sec:newidQM}.
The potential function $V_2(x)$ can be obtained from \eqref{V=V1V2}.

\noindent
Potential functions $V_1(x;\bm{\lambda})$:
\begin{align}
%%%%%%%%%%%%%%%%%%%%%%%%%%%%%%%%
% continuous Hahn (cH)         %
%%%%%%%%%%%%%%%%%%%%%%%%%%%%%%%%
  \text{cH}&:\text{(a)}:
  \ V_1(x;\bm{\lambda})=a_1+ix,\quad
  \text{(b)}:
  \ V_1(x;\bm{\lambda})=a_2+ix,\\
%
%%%%%%%%%%%%%%%%%%%%%%%%%%%%%%%%
% Meixner-Pollaczek (MP)       %
%%%%%%%%%%%%%%%%%%%%%%%%%%%%%%%%
  \text{MP}&:\text{(a)}:
  \ V_1(x;\bm{\lambda})=a+ix,\quad
  \text{(b)}:
  \ V_1(x;\bm{\lambda})=e^{i(\frac{\pi}{2}-\phi)},\\
%
%%%%%%%%%%%%%%%%%%%%%%%%%%%%%%%%
% Wilson (W)                   %
%%%%%%%%%%%%%%%%%%%%%%%%%%%%%%%%
  \text{W}&:
  \ \text{Assume $\{a_j^*,a_k^*\}=\{a_j,a_k\}$ (as a set) and
  set $\{l,m\}=\{1,2,3,4\}\backslash\{j,k\}$,}\n
  &\ \ \ V_i(x;\bm{\lambda})=V^{(j,k)}_i(x;\bm{\lambda})\ \ (i=1,2),\quad
  V_1(x;\bm{\lambda})=\frac{(a_j+ix)(a_k+ix)}{2ix+1},\\
%
%%%%%%%%%%%%%%%%%%%%%%%%%%%%%%%%
% continuous dual Hahn (cdH)   %
%%%%%%%%%%%%%%%%%%%%%%%%%%%%%%%%
  \text{cdH}&:
  \ \text{Assume $\{a_j^*,a_k^*\}=\{a_j,a_k\}$ (as a set) and
  set $\{l\}=\{1,2,3\}\backslash\{j,k\}$,}\n
  &\ \ \ V_i(x;\bm{\lambda})=V^{(j,k)}_i(x;\bm{\lambda})\ \ (i=1,2),\n
  &\ \ \ \text{(a)}:
  \ V_1(x;\bm{\lambda})=\frac{(a_j+ix)(a_k+ix)}{2ix+1},\quad
  \text{(b)}:
  \ V_1(x;\bm{\lambda})=\frac{a_l+ix}{2ix+1},\\
%
%%%%%%%%%%%%%%%%%%%%%%%%%%%%%%%%%%%%%
% continuous dual $q$-Hahn (cd$q$H) %
%%%%%%%%%%%%%%%%%%%%%%%%%%%%%%%%%%%%%
  \text{cd$q$H}&:
  \ \text{Assume $\{a_j^*,a_k^*\}=\{a_j,a_k\}$ (as a set) and
  set $\{l\}=\{1,2,3\}\backslash\{j,k\}$,}\n
  &\ \ \ V_i(x;\bm{\lambda})=V^{(j,k)}_i(x;\bm{\lambda})\ \ (i=1,2),\n
  &\ \ \ \text{(a)}:
  \ V_1(x;\bm{\lambda})=\frac{(1-a_je^{ix})(1-a_ke^{ix})}{1-qe^{2ix}},\quad
  \text{(b)}:
  \ V_1(x;\bm{\lambda})=\frac{1-a_le^{ix}}{1-qe^{2ix}},\\
%
%%%%%%%%%%%%%%%%%%%%%%%%%%%%%%%%
% Al-Salam-Chihara (ASC)       %
%%%%%%%%%%%%%%%%%%%%%%%%%%%%%%%%
  \text{ASC}&:
  \ \text{Assume $a_1,a_2\in\mathbb{R}$ for (b) and (c),}\n
  &\ \ \ \text{(a)}:
  \ V_1(x;\bm{\lambda})=\frac{(1-a_1e^{ix})(1-a_2e^{ix})}{1-qe^{2ix}},\quad
  \text{(b)}:
  \ V_1(x;\bm{\lambda})=\frac{1-a_1e^{ix}}{1-qe^{2ix}},\\
  &\ \ \ \text{(c)}:
  \ V_1(x;\bm{\lambda})=\frac{1-a_2e^{ix}}{1-qe^{2ix}},\quad
  \text{(d)}:
  \ V_1(x;\bm{\lambda})=\frac{1}{1-qe^{2ix}},\\
%
%%%%%%%%%%%%%%%%%%%%%%%%%%%%%%%%%%%%%%%%
% continuous big $q$-Hermite (cb$q$He) %
%%%%%%%%%%%%%%%%%%%%%%%%%%%%%%%%%%%%%%%%
  \text{cb$q$He}&:\text{(a)}:
  \ V_1(x;\bm{\lambda})=\frac{1-ae^{ix}}{1-qe^{2ix}},\quad
  \text{(b)}:
  \ V_1(x;\bm{\lambda})=\frac{1}{1-qe^{2ix}},\\
%
%%%%%%%%%%%%%%%%%%%%%%%%%%%%%%%%%%%
% continuous $q$-Hermite (c$q$He) %
%%%%%%%%%%%%%%%%%%%%%%%%%%%%%%%%%%%
  \text{c$q$He}&:
  \ V_1(x;\bm{\lambda})=\frac{1}{1-qe^{2ix}},\\
%
%%%%%%%%%%%%%%%%%%%%%%%%%%%%%%%%%
% continuous $q$-Jacobi (c$q$J) %
%%%%%%%%%%%%%%%%%%%%%%%%%%%%%%%%%
  \text{c$q$J}&:\text{(a)}:
  \ V_1(x;\bm{\lambda})=\frac{(1-q^{\frac12(\alpha+\frac12)}e^{ix})
  (1-q^{\frac12(\alpha+\frac32)}e^{ix})}{1-qe^{2ix}},\\
  &\ \ \,\text{(b)}:
  \ V_1(x;\bm{\lambda})=\frac{(1+q^{\frac12(\beta+\frac12)}e^{ix})
  (1+q^{\frac12(\beta+\frac32)}e^{ix})}{1-qe^{2ix}},\\
%
%%%%%%%%%%%%%%%%%%%%%%%%%%%%%%%%%%%
% continuous $q$-Laguerre (c$q$L) %
%%%%%%%%%%%%%%%%%%%%%%%%%%%%%%%%%%%
  \text{c$q$L}&:\text{(a)}:
  \ V_1(x;\bm{\lambda})=\frac{(1-q^{\frac12(\alpha+\frac12)}e^{ix})
  (1-q^{\frac12(\alpha+\frac32)}e^{ix})}{1-qe^{2ix}},
  \ \ \text{(b)}:
  \ V_1(x;\bm{\lambda})=\frac{1}{1-qe^{2ix}},\\
%
%%%%%%%%%%%%%%%%%%%%%%%%%%%%%%%
% continuous $q$-Hahn (c$q$H) %
%%%%%%%%%%%%%%%%%%%%%%%%%%%%%%%
  \text{c$q$H}&:\text{(a)}:
  \ V_1(x;\bm{\lambda})=\frac{(1-a_1e^{2i\phi}e^{ix})(1-a_1^*e^{ix})}
  {1-qe^{2i\phi}e^{2ix}},\\
  &\ \ \,\text{(b)}:
  \ V_1(x;\bm{\lambda})=\frac{(1-a_2e^{2i\phi}e^{ix})(1-a_2^*e^{ix})}
  {1-qe^{2i\phi}e^{2ix}},\\
%
%%%%%%%%%%%%%%%%%%%%%%%%%%%%%%%%%
% $q$-Meixner-Pollaczek ($q$MP) %
%%%%%%%%%%%%%%%%%%%%%%%%%%%%%%%%%
  \text{$q$MP}&:\text{(a)}:
  \ V_1(x;\bm{\lambda})=\frac{(1-ae^{2i\phi}e^{ix})(1-ae^{ix})}
  {1-qe^{2i\phi}e^{2ix}},\quad
  \text{(b)}:
  \ V_1(x;\bm{\lambda})=\frac{1}{1-qe^{2i\phi}e^{2ix}}.
\end{align}
Constants $\tilde{f}_n(\bm{\lambda})$, $\tilde{b}_n(\bm{\lambda})$ and
$\bm{\bar{\delta}}$:
\begin{align}
%%%%%%%%%%%%%%%%%%%%%%%%%%%%%%%%
% continuous Hahn (cH)         %
%%%%%%%%%%%%%%%%%%%%%%%%%%%%%%%%
  \text{cH}&:\text{(a)}:
  \ \tilde{f}_n(\bm{\lambda})=a_1+a_1^*+n-1,
  \ \ \tilde{b}_n(\bm{\lambda})=a_2+a_2^*+n,
  \ \ \bm{\bar{\delta}}=(\tfrac12,-\tfrac12),\\
  &\ \ \,\text{(b)}:
  \ \tilde{f}_n(\bm{\lambda})=a_2+a_2^*+n-1,
  \ \ \tilde{b}_n(\bm{\lambda})=a_1+a_1^*+n,
  \ \ \bm{\bar{\delta}}=(-\tfrac12,\tfrac12),\\
%
%%%%%%%%%%%%%%%%%%%%%%%%%%%%%%%%
% Meixner-Pollaczek (MP)       %
%%%%%%%%%%%%%%%%%%%%%%%%%%%%%%%%
  \text{MP}&:\text{(a)}:
  \ \tilde{f}_n(\bm{\lambda})=2a+n-1,
  \ \ \tilde{b}_n(\bm{\lambda})=2\sin\phi,
  \ \ \bm{\bar{\delta}}=(\tfrac12,0),\\
  &\ \ \,\text{(b)}:
  \ \tilde{f}_n(\bm{\lambda})=2\sin\phi,
  \ \ \tilde{b}_n(\bm{\lambda})=2a+n,
  \ \ \bm{\bar{\delta}}=(-\tfrac12,0),\\
%
%%%%%%%%%%%%%%%%%%%%%%%%%%%%%%%%
% Wilson (W)                   %
%%%%%%%%%%%%%%%%%%%%%%%%%%%%%%%%
  \text{W}&:
  \ \tilde{f}_n(\bm{\lambda})=a_j+a_k+n-1,
  \ \ \tilde{b}_n(\bm{\lambda})=a_l+a_m+n,\n
  &\ \ \ (\bm{\bar{\delta}})_j=(\bm{\bar{\delta}})_k=\tfrac12,
  \ \ (\bm{\bar{\delta}})_l=(\bm{\bar{\delta}})_m=-\tfrac12,\\
%
%%%%%%%%%%%%%%%%%%%%%%%%%%%%%%%%
% continuous dual Hahn (cdH)   %
%%%%%%%%%%%%%%%%%%%%%%%%%%%%%%%%
  \text{cdH}&:\text{(a)}:
  \ \tilde{f}_n(\bm{\lambda})=a_j+a_k+n-1,
  \ \ \tilde{b}_n(\bm{\lambda})=1,
  \ \ (\bm{\bar{\delta}})_j=(\bm{\bar{\delta}})_k=\tfrac12,
  \ \ (\bm{\bar{\delta}})_l=-\tfrac12,\\
  &\ \ \,\text{(b)}:
  \ \tilde{f}_n(\bm{\lambda})=1,
  \ \ \tilde{b}_n(\bm{\lambda})=a_j+a_k+n,
  \ \ (\bm{\bar{\delta}})_l=\tfrac12,
  \ \ (\bm{\bar{\delta}})_j=(\bm{\bar{\delta}})_k=-\tfrac12,\\
%
%%%%%%%%%%%%%%%%%%%%%%%%%%%%%%%%%%%%%
% continuous dual $q$-Hahn (cd$q$H) %
%%%%%%%%%%%%%%%%%%%%%%%%%%%%%%%%%%%%%
  \text{cd$q$H}&:\text{(a)}:
  \tilde{f}_n(\bm{\lambda})=q^{-\frac{n}{2}}(1-a_ja_kq^{n-1}),
  \ \tilde{b}_n(\bm{\lambda})=q^{-\frac{n}{2}},
  \,(\bm{\bar{\delta}})_j=(\bm{\bar{\delta}})_k=\tfrac12,
  \,(\bm{\bar{\delta}})_l=-\tfrac12,\!\\
  &\ \ \,\text{(b)}:
  \tilde{f}_n(\bm{\lambda})=q^{-\frac{n}{2}},
  \ \tilde{b}_n(\bm{\lambda})=q^{-\frac{n}{2}}(1-a_ja_kq^n),
  \ (\bm{\bar{\delta}})_l=\tfrac12,
  \ (\bm{\bar{\delta}})_j=(\bm{\bar{\delta}})_k=-\tfrac12,\\
%
%%%%%%%%%%%%%%%%%%%%%%%%%%%%%%%%
% Al-Salam-Chihara (ASC)       %
%%%%%%%%%%%%%%%%%%%%%%%%%%%%%%%%
  \text{ASC}&:\text{(a)}:
  \ \tilde{f}_n(\bm{\lambda})=q^{-\frac{n}{2}}(1-a_1a_2q^{n-1}),
  \ \ \tilde{b}_n(\bm{\lambda})=q^{-\frac{n}{2}},
  \ \ \bm{\bar{\delta}}=(\tfrac12,\tfrac12),\\
  &\ \ \,\text{(b)}:
  \ \tilde{f}_n(\bm{\lambda})=q^{-\frac{n}{2}},
  \ \ \tilde{b}_n(\bm{\lambda})=q^{-\frac{n}{2}},
  \ \ \bm{\bar{\delta}}=(\tfrac12,-\tfrac12),\\
  &\ \ \,\text{(c)}:
  \ \tilde{f}_n(\bm{\lambda})=q^{-\frac{n}{2}},
  \ \ \tilde{b}_n(\bm{\lambda})=q^{-\frac{n}{2}},
  \ \ \bm{\bar{\delta}}=(-\tfrac12,\tfrac12),\\
  &\ \ \,\text{(d)}:
  \ \tilde{f}_n(\bm{\lambda})=q^{-\frac{n}{2}},
  \ \ \tilde{b}_n(\bm{\lambda})=q^{-\frac{n}{2}}(1-a_1a_2q^n),
  \ \ \bm{\bar{\delta}}=(-\tfrac12,-\tfrac12),\\
%
%%%%%%%%%%%%%%%%%%%%%%%%%%%%%%%%%%%%%%%%
% continuous big $q$-Hermite (cb$q$He) %
%%%%%%%%%%%%%%%%%%%%%%%%%%%%%%%%%%%%%%%%
  \text{cb$q$He}&:\text{(a)}:
  \ \tilde{f}_n(\bm{\lambda})=q^{-\frac{n}{2}},
  \ \ \tilde{b}_n(\bm{\lambda})=q^{-\frac{n}{2}},
  \ \ \bm{\bar{\delta}}=\tfrac12,\\
  &\ \ \,\text{(b)}:
  \ \tilde{f}_n(\bm{\lambda})=q^{-\frac{n}{2}},
  \ \ \tilde{b}_n(\bm{\lambda})=q^{-\frac{n}{2}},
  \ \ \bm{\bar{\delta}}=-\tfrac12,\\
%
%%%%%%%%%%%%%%%%%%%%%%%%%%%%%%%%%%%
% continuous $q$-Hermite (c$q$He) %
%%%%%%%%%%%%%%%%%%%%%%%%%%%%%%%%%%%
  \text{c$q$He}&:
  \ \tilde{f}_n(\bm{\lambda})=q^{-\frac{n}{2}},
  \ \ \tilde{b}_n(\bm{\lambda})=q^{-\frac{n}{2}},
  \ \ \bm{\bar{\delta}}:\text{none},\\
%
%%%%%%%%%%%%%%%%%%%%%%%%%%%%%%%%%
% continuous $q$-Jacobi (c$q$J) %
%%%%%%%%%%%%%%%%%%%%%%%%%%%%%%%%%
  \text{c$q$J}&:\text{(a)}:
  \ \tilde{f}_n(\bm{\lambda})=1-q^{\alpha+n},
  \ \ \tilde{b}_n(\bm{\lambda})=q^{-n}(1-q^{\beta+n+1}),
  \ \ \bm{\bar{\delta}}=(1,-1),\\
  &\ \ \,\text{(b)}:
  \ \tilde{f}_n(\bm{\lambda})=q^{-n}(1-q^{\beta+n}),
  \ \ \tilde{b}_n(\bm{\lambda})=1-q^{\alpha+n+1},
  \ \ \bm{\bar{\delta}}=(-1,1),\\
%
%%%%%%%%%%%%%%%%%%%%%%%%%%%%%%%%%%%
% continuous $q$-Laguerre (c$q$L) %
%%%%%%%%%%%%%%%%%%%%%%%%%%%%%%%%%%%
  \text{c$q$L}&:\text{(a)}:
  \ \tilde{f}_n(\bm{\lambda})=1-q^{\alpha+n},
  \ \ \tilde{b}_n(\bm{\lambda})=q^{-n},
  \ \ \bm{\bar{\delta}}=1,\\
  &\ \ \,\text{(b)}:
  \ \tilde{f}_n(\bm{\lambda})=q^{-n},
  \ \ \tilde{b}_n(\bm{\lambda})=1-q^{\alpha+n+1},
  \ \ \bm{\bar{\delta}}=-1,\\
%
%%%%%%%%%%%%%%%%%%%%%%%%%%%%%%%
% continuous $q$-Hahn (c$q$H) %
%%%%%%%%%%%%%%%%%%%%%%%%%%%%%%%
  \text{c$q$H}&:\text{(a)}:
  \ \tilde{f}_n(\bm{\lambda})=q^{-\frac{n}{2}}(1-a_1a_1^*q^{n-1}),
  \,\tilde{b}_n(\bm{\lambda})=q^{-\frac{n}{2}}(1-a_2a_2^*q^n),
  \,\bm{\bar{\delta}}=(\tfrac12,-\tfrac12,0),\\
  &\ \ \,\text{(b)}:
  \ \tilde{f}_n(\bm{\lambda})=q^{-\frac{n}{2}}(1-a_2a_2^*q^{n-1}),
  \,\tilde{b}_n(\bm{\lambda})=q^{-\frac{n}{2}}(1-a_1a_1^*q^n),
  \,\bm{\bar{\delta}}=(-\tfrac12,\tfrac12,0),\\
%
%%%%%%%%%%%%%%%%%%%%%%%%%%%%%%%%%
% $q$-Meixner-Pollaczek ($q$MP) %
%%%%%%%%%%%%%%%%%%%%%%%%%%%%%%%%%
  \text{$q$MP}&:\text{(a)}:
  \ \tilde{f}_n(\bm{\lambda})=q^{-\frac{n}{2}}(1-a^2q^{n-1}),
  \ \ \tilde{b}_n(\bm{\lambda})=q^{-\frac{n}{2}},
  \ \ \bm{\bar{\delta}}=(\tfrac12,0),\\
  &\ \ \,\text{(b)}:
  \ \tilde{f}_n(\bm{\lambda})=q^{-\frac{n}{2}},
  \ \ \tilde{b}_n(\bm{\lambda})=q^{-\frac{n}{2}}(1-a^2q^n),
  \ \ \bm{\bar{\delta}}=(-\tfrac12,0).
\end{align}

%%%%%%%%%%%%%%%%%%%%%%%%%%%%%%%%%%%%%%%%%%%%%%%%%%%%%
%                                                   %
% A.3 Data for \S\,\ref{sec:newrdQM}                %
%                                                   %
%%%%%%%%%%%%%%%%%%%%%%%%%%%%%%%%%%%%%%%%%%%%%%%%%%%%%
\subsection{Data for \S\,\ref{sec:newrdQM}}
\label{app:rdQM}

We present explicit forms of $B_1(x)$, $D_1(x)$, $\tilde{f}_n$, $\tilde{b}_n$
and $\bm{\bar{\delta}}$ in \S\,\ref{sec:newrdQM}.
The potential functions $B_2(x)$ and $D_2(x)$ can be obtained from
\eqref{B=B1B2}.

\noindent
Potential functions $B_1(x;\bm{\lambda})$ and $D_1(x;\bm{\lambda})$:
\begin{align}
%%%%%%%%%%%%%%%%%%%%%%%%%%%%%%%%
% Hahn (H)                     %
%%%%%%%%%%%%%%%%%%%%%%%%%%%%%%%%
  \text{H}&:\text{(a)}:
  \ B_1(x;\bm{\lambda})=\frac{N-x}{N+1},
  \ \ D_1(x;\bm{\lambda})=\frac{x}{N+1},\\
  &\ \ \,\text{(b)}:
  \ B_1(x;\bm{\lambda})=\frac{x+a}{a-1},
  \ \ D_1(x;\bm{\lambda})=\frac{x}{1-a},\\
  &\ \ \,\text{(c)}:
  \ B_1(x;\bm{\lambda})=a(N-x),
  \ \ D_1(x;\bm{\lambda})=-a(b+N-x),\\
  &\ \ \,\text{(d)}:
  \ B_1(x;\bm{\lambda})=N(x+a),
  \ \ D_1(x;\bm{\lambda})=N(b+N-x),\\
%
%%%%%%%%%%%%%%%%%%%%%%%%%%%%%%%%
% Krawtchouk (K)               %
%%%%%%%%%%%%%%%%%%%%%%%%%%%%%%%%
  \text{K}&:\text{(a)}:
  \ B_1(x;\bm{\lambda})=\frac{N-x}{N+1},
  \ \ D_1(x;\bm{\lambda})=\frac{x}{N+1},\\
  &\ \ \,\text{(b)}:
  \ B_1(x;\bm{\lambda})=Np,
  \ \ D_1(x;\bm{\lambda})=N(1-p),\\
%
%%%%%%%%%%%%%%%%%%%%%%%%%%%%%%%%
% Racah (R)                    %
%%%%%%%%%%%%%%%%%%%%%%%%%%%%%%%%
  \text{R}&:\text{(a)}:
  \ B_1(x;\bm{\lambda})=-\frac{(x-N)(x+d)}{(N+1)(2x+1+d)},
  \ \ D_1(x;\bm{\lambda})=\frac{(x+d+N)x}{(N+1)(2x-1+d)},\!\\
  &\ \ \,\text{(b)}:
  \ B_1(x;\bm{\lambda})=\frac{(x+b)(x+d)}{(b-1)(2x+1+d)},
  \ \ D_1(x;\bm{\lambda})=-\frac{(x+d-b)x}{(b-1)(2x-1+d)},\\
  &\ \ \,\text{(c)}:
  \ B_1(x;\bm{\lambda})=\frac{(x+c)(x+d)}{(c-1)(2x+1+d)},
  \ \ D_1(x;\bm{\lambda})=-\frac{(x+d-c)x}{(c-1)(2x-1+d)},\\
  &\ \ \,\text{(d)}:
  \ B_1(x;\bm{\lambda})=-\frac{c(x-N)(x+b)}{2x+1+d},
  \ \ D_1(x;\bm{\lambda})=\frac{c(x+d+N)(x+d-b)}{2x-1+d},\\
  &\ \ \,\text{(e)}:
  \ B_1(x;\bm{\lambda})=-\frac{b(x-N)(x+c)}{2x+1+d},
  \ \ D_1(x;\bm{\lambda})=\frac{b(x+d+N)(x+d-c)}{2x-1+d},\\
  &\ \ \,\text{(f)}:
  \ B_1(x;\bm{\lambda})=\frac{N(x+b)(x+c)}{2x+1+d},
  \ \ D_1(x;\bm{\lambda})=-\frac{N(x+d-b)(x+d-c)}{2x-1+d},\\
%
%%%%%%%%%%%%%%%%%%%%%%%%%%%%%%%%
% dual Hahn (dH)               %
%%%%%%%%%%%%%%%%%%%%%%%%%%%%%%%%
  \text{dH}&:\text{(a)}:
  \ B_1(x;\bm{\lambda})=\frac{(x+a+b-1)(N-x)}{(N+1)(2x+a+b)},
  \ \ D_1(x;\bm{\lambda})=\frac{x(x+a+b+N-1)}{(N+1)(2x-2+a+b)},\\
  &\ \ \,\text{(b)}:
  \ B_1(x;\bm{\lambda})=\frac{(x+a)(x+a+b-1)}{(a-1)(2x+a+b)},
  \ \ D_1(x;\bm{\lambda})=\frac{x(x+b-1)}{(1-a)(2x-2+a+b)},\\
  &\ \ \,\text{(c)}:
  \ B_1(x;\bm{\lambda})=\frac{x+a+b-1}{2x+a+b},
  \ \ D_1(x;\bm{\lambda})=\frac{x}{2x-2+a+b},\\
  &\ \ \,\text{(d)}:
  \ B_1(x;\bm{\lambda})=\frac{(x+a)(N-x)}{2x+a+b},
  \ \ D_1(x;\bm{\lambda})=\frac{(x+b-1)(x+a+b+N-1)}{2x-2+a+b},\\
  &\ \ \,\text{(e)}:
  \ B_1(x;\bm{\lambda})=\frac{a(N-x)}{2x+a+b},
  \ \ D_1(x;\bm{\lambda})=-\frac{a(x+a+b+N-1)}{2x-2+a+b},\\
  &\ \ \,\text{(f)}:
  \ B_1(x;\bm{\lambda})=\frac{N(x+a)}{2x+a+b},
  \ \ D_1(x;\bm{\lambda})=\frac{N(x+b-1)}{2x-2+a+b},\\
%
%%%%%%%%%%%%%%%%%%%%%%%%%%%%%%%%%%%%%%
% dual quantum q-Krawtchouk (dq$q$K) %
%%%%%%%%%%%%%%%%%%%%%%%%%%%%%%%%%%%%%%
  \text{dq$q$K}&:\text{(a)}:
  \ B_1(x;\bm{\lambda})=\frac{q^{-N-1}(1-q^{N-x})}{q^{-N-1}-1},
  \ \ D_1(x;\bm{\lambda})=\frac{q^{-x}-1}{q^{-N-1}-1},\\
  &\ \ \,\text{(b)}:
  \ B_1(x;\bm{\lambda})=q^{-x-1},
  \ \ D_1(x;\bm{\lambda})=-(q^{-x}-1),\\
  &\ \ \,\text{(c)}:
  \ B_1(x;\bm{\lambda})=p^{-1}q^{-N-1}(1-q^{N-x}),
  \ \ D_1(x;\bm{\lambda})=-(1-p^{-1}q^{-x}),\\
  &\ \ \,\text{(d)}:
  \ B_1(x;\bm{\lambda})=(1-q^N)p^{-1}q^{-x-N-1},
  \ \ D_1(x;\bm{\lambda})=(q^{-N}-1)(1-p^{-1}q^{-x}),\\
%
%%%%%%%%%%%%%%%%%%%%%%%%%%%%%%%%
% q-Hahn (qH)                  %
%%%%%%%%%%%%%%%%%%%%%%%%%%%%%%%%
  \text{$q$H}&:\text{(a)}:
  \ B_1(x;\bm{\lambda})=\frac{q^{x-N}-1}{q^{-N-1}-1},
  \ \ D_1(x;\bm{\lambda})=\frac{1-q^x}{1-q^{N+1}},\\
  &\ \ \,\text{(b)}:
  \ B_1(x;\bm{\lambda})=\frac{1-aq^x}{1-aq^{-1}},
  \ \ D_1(x;\bm{\lambda})=\frac{aq^{-1}(1-q^x)}{aq^{-1}-1},\\
  &\ \ \,\text{(c)}:
  \ B_1(x;\bm{\lambda})=(1-a)(q^{x-N}-1),
  \ \ D_1(x;\bm{\lambda})=(a-1)q^{-1}(q^{x-N}-b),\\
  &\ \ \,\text{(d)}:
  \ B_1(x;\bm{\lambda})=(q^{-N}-1)(1-aq^x),
  \ \ D_1(x;\bm{\lambda})=(1-q^N)aq^{-1}(q^{x-N}-b),\\
%
%%%%%%%%%%%%%%%%%%%%%%%%%%%%%%%%
% q-Krawtchouk ($q$K)          %
%%%%%%%%%%%%%%%%%%%%%%%%%%%%%%%%
  \text{$q$K}&:\text{(a)}:
  \ B_1(x;\bm{\lambda})=\frac{q^{x-N}-1}{q^{-N-1}-1},
  \ \ D_1(x;\bm{\lambda})=\frac{1-q^x}{1-q^{N+1}},\\
  &\ \ \,\text{(b)}:
  \ B_1(x;\bm{\lambda})=q^{-N}-1,
  \ \ D_1(x;\bm{\lambda})=(1-q^N)p,\\
%
%%%%%%%%%%%%%%%%%%%%%%%%%%%%%%%%
% quantum q-Krawtchouk (q$q$K) %
%%%%%%%%%%%%%%%%%%%%%%%%%%%%%%%%
  \text{q$q$K}&:\text{(a)}:
  \ B_1(x;\bm{\lambda})=\frac{q^{x-N}-1}{q^{-N-1}-1},
  \ \ D_1(x;\bm{\lambda})=\frac{1-q^x}{1-q^{N+1}},\\
  &\ \ \,\text{(b)}:
  \ B_1(x;\bm{\lambda})=qq^x,
  \ \ D_1(x;\bm{\lambda})=1-q^x,\\
  &\ \ \,\text{(c)}:
  \ B_1(x;\bm{\lambda})=p^{-1}(q^{x-N}-1),
  \ \ D_1(x;\bm{\lambda})=1-p^{-1}q^{x-N-1},\\
  &\ \ \,\text{(d)}:
  \ B_1(x;\bm{\lambda})=(q^{-N}-1)p^{-1}q^x,
  \ \ D_1(x;\bm{\lambda})=(1-q^N)(1-p^{-1}q^{x-N-1}),\\
%
%%%%%%%%%%%%%%%%%%%%%%%%%%%%%%%%
% affine q-Krawtchouk (a$q$K)  %
%%%%%%%%%%%%%%%%%%%%%%%%%%%%%%%%
  \text{a$q$K}&:\text{(a)}:
  \ B_1(x;\bm{\lambda})=\frac{q^{x-N}-1}{q^{-N-1}-1},
  \ \ D_1(x;\bm{\lambda})=\frac{1-q^x}{1-q^{N+1}},\\
  &\ \ \,\text{(b)}:
  \ B_1(x;\bm{\lambda})=\frac{1-pq^{x+1}}{1-p},
  \ \ D_1(x;\bm{\lambda})=\frac{p(1-q^x)}{p-1},\\
  &\ \ \,\text{(c)}:
  \ B_1(x;\bm{\lambda})=(1-pq)(q^{x-N}-1),
  \ \ D_1(x;\bm{\lambda})=(pq-1)q^{x-N-1},\\
  &\ \ \,\text{(d)}:
  \ B_1(x;\bm{\lambda})=(q^{-N}-1)(1-pq^{x+1}),
  \ \ D_1(x;\bm{\lambda})=(q^{-N}-1)pq^x,\\
%
%%%%%%%%%%%%%%%%%%%%%%%%%%%%%%%%
% dual q-Hahn (dqH)            %
%%%%%%%%%%%%%%%%%%%%%%%%%%%%%%%%
  \text{d$q$H}&:\text{(a)}:
  \ B_1(x;\bm{\lambda})=\frac{(q^{x-N}-1)(1-abq^{x-1})}
  {(q^{-N-1}-1)(1-abq^{2x})},
  \ \ D_1(x;\bm{\lambda})=\frac{(1-q^x)(1-abq^{x+N-1})}
  {(1-q^{N+1})(1-abq^{2x-2})},\\
  &\ \ \,\text{(b)}:
  \ B_1(x;\bm{\lambda})=\frac{(1-aq^x)(1-abq^{x-1})}{(1-aq^{-1})(1-abq^{2x})},
  \ \ D_1(x;\bm{\lambda})=\frac{(1-q^x)(1-bq^{x-1})}
  {(1-a^{-1}q)(1-abq^{2x-2})},\\
  &\ \ \,\text{(c)}:
  \ B_1(x;\bm{\lambda})=\frac{1-abq^{x-1}}{1-abq^{2x}},
  \ \ D_1(x;\bm{\lambda})=\frac{bq^{-2}aq^x(1-q^x)}{1-abq^{2x-2}},\\
  &\ \ \,\text{(d)}:
  \ B_1(x;\bm{\lambda})=\frac{(q^{x-N}-1)(1-aq^x)}{1-abq^{2x}},
  \ \ D_1(x;\bm{\lambda})=\frac{q^{-N}(1-abq^{x+N-1})(1-bq^{x-1})}
  {b(1-abq^{2x-2})},\\
  &\ \ \,\text{(e)}:
  \ B_1(x;\bm{\lambda})=\frac{(1-a)(q^{x-N}-1)}{1-abq^{2x}},
  \ \ D_1(x;\bm{\lambda})=\frac{(a-1)q^{x-N-1}(1-abq^{x+N-1})}{1-abq^{2x-2}},\\
  &\ \ \,\text{(f)}:
  \ B_1(x;\bm{\lambda})=\frac{(q^{-N}-1)(1-aq^x)}{1-abq^{2x}},
  \ \ D_1(x;\bm{\lambda})=\frac{(1-q^N)aq^{x-N-1}(1-bq^{x-1})}{1-abq^{2x-2}},\\
%
%%%%%%%%%%%%%%%%%%%%%%%%%%%%%%%%
% dual q-Krawtchouk (dqK)      %
%%%%%%%%%%%%%%%%%%%%%%%%%%%%%%%%
  \text{d$q$K}&:\text{(a)}:
  \ B_1(x;\bm{\lambda})=\frac{(q^{x-N}-1)(1+pq^x)}{(q^{-N-1}-1)(1+pq^{2x+1})},
  \ \ D_1(x;\bm{\lambda})=\frac{(1-q^x)(1+pq^{x+N})}{(1-q^{N+1})(1+pq^{2x-1})},\\
  &\ \ \,\text{(b)}:
  \ B_1(x;\bm{\lambda})=\frac{1+pq^x}{1+pq^{2x+1}},
  \ \ D_1(x;\bm{\lambda})=-\frac{pq^{x-1}(1-q^x)}{1+pq^{2x-1}},\\
  &\ \ \,\text{(c)}:
  \ B_1(x;\bm{\lambda})=\frac{q^{x-N}-1}{1+pq^{2x+1}},
  \ \ D_1(x;\bm{\lambda})=-\frac{q^{x-N-1}(1+pq^{x+N})}{1+pq^{2x-1}},\\
  &\ \ \,\text{(d)}:
  \ B_1(x;\bm{\lambda})=\frac{q^{-N}-1}{1+pq^{2x+1}},
  \ \ D_1(x;\bm{\lambda})=\frac{(1-q^N)pq^{2x-N-1}}{1+pq^{2x-1}},\\
%
%%%%%%%%%%%%%%%%%%%%%%%%%%%%%%%%
% Meixner (M)                  %
%%%%%%%%%%%%%%%%%%%%%%%%%%%%%%%%
  \text{M}&:\text{(a)}:
  \ B_1(x;\bm{\lambda})=\frac{x+\beta}{\beta-1},
  \ \ D_1(x;\bm{\lambda})=\frac{x}{1-\beta},\\
  &\ \ \,\text{(b)}:
  \ B_1(x;\bm{\lambda})=\frac{\beta c}{1-c},
  \ \ D_1(x;\bm{\lambda})=-\frac{\beta}{1-c},\\
%
%%%%%%%%%%%%%%%%%%%%%%%%%%%%%%%%
% little $q$-Jacobi (l$q$J)    %
%%%%%%%%%%%%%%%%%%%%%%%%%%%%%%%%
  \text{l$q$J}&:\text{(a)}:
  \ B_1(x;\bm{\lambda})=\frac{q^{-1}(q^{-x}-b)}{1-bq^{-1}},
  \ \ D_1(x;\bm{\lambda})=\frac{q^{-x}-1}{bq^{-1}-1},\\
  &\ \ \,\text{(b)}:
  \ B_1(x;\bm{\lambda})=(1-b)aq^{-1},
  \ \ D_1(x;\bm{\lambda})=b-1,\\
%
%%%%%%%%%%%%%%%%%%%%%%%%%%%%%%%%
% little $q$-Laguerre (l$q$L)  %
%%%%%%%%%%%%%%%%%%%%%%%%%%%%%%%%
  \text{l$q$L}&:\text{(a)}:
  \ B_1(x;\bm{\lambda})=q^{-x-1},
  \ \ D_1(x;\bm{\lambda})=-(q^{-x}-1),\\
  &\ \ \,\text{(b)}:
  \ B_1(x;\bm{\lambda})=aq^{-1},
  \ \ D_1(x;\bm{\lambda})=-1,\\
%
%%%%%%%%%%%%%%%%%%%%%%%%%%%%%%%%
% $q$-Meixner ($q$M)           %
%%%%%%%%%%%%%%%%%%%%%%%%%%%%%%%%
  \text{$q$M}&:\text{(a)}:
  \ B_1(x;\bm{\lambda})=q^{x+1},
  \ \ D_1(x;\bm{\lambda})=1-q^x,\\
  &\ \ \,\text{(b)}:
  \ B_1(x;\bm{\lambda})=\frac{1-bq^{x+1}}{1-b},
  \ \ D_1(x;\bm{\lambda})=\frac{1-q^x}{1-b^{-1}},\\
  &\ \ \,\text{(c)}:
  \ B_1(x;\bm{\lambda})=(1-bq)cq^x,
  \ \ D_1(x;\bm{\lambda})=(1-b^{-1}q^{-1})(1+bcq^x),\\
  &\ \ \,\text{(d)}:
  \ B_1(x;\bm{\lambda})=c(1-bq^{x+1}),
  \ \ D_1(x;\bm{\lambda})=1+bcq^x,\\
%
%%%%%%%%%%%%%%%%%%%%%%%%%%%%%%%%%%%%%
% Al-Salam-Carlitz $\II$ (ASC$\II$) %
%%%%%%%%%%%%%%%%%%%%%%%%%%%%%%%%%%%%%
  \text{ASC$\II$}&:\text{(a)}:
  \ B_1(x;\bm{\lambda})=q^{x+1},
  \ \ D_1(x;\bm{\lambda})=1-q^x,\\
  &\ \ \,\text{(b)}:
  \ B_1(x;\bm{\lambda})=aq^{x+1},
  \ \ D_1(x;\bm{\lambda})=1-aq^x,\\
%
%%%%%%%%%%%%%%%%%%%%%%%%%%%%%%%%
% $q$-Charlier ($q$C)          %
%%%%%%%%%%%%%%%%%%%%%%%%%%%%%%%%
  \text{$q$C}&:\text{(a)}:
  \ B_1(x;\bm{\lambda})=q^{x+1},
  \ \ D_1(x;\bm{\lambda})=1-q^x,\\
  &\ \ \,\text{(b)}:
  \ B_1(x;\bm{\lambda})=a,
  \ \ D_1(x;\bm{\lambda})=1.
\end{align}
Constants $\tilde{f}_n(\bm{\lambda})$, $\tilde{b}_n(\bm{\lambda})$ and
$\bm{\bar{\delta}}$:
\begin{align}
%%%%%%%%%%%%%%%%%%%%%%%%%%%%%%%%
% Hahn (H)                     %
%%%%%%%%%%%%%%%%%%%%%%%%%%%%%%%%
  \text{H}&:\text{(a)}:
  \ \tilde{f}_n(\bm{\lambda})=1,
  \ \ \tilde{b}_n(\bm{\lambda})
  =\mathcal{E}_{N+1}(\bm{\lambda})-\mathcal{E}_n(\bm{\lambda}),
  \ \ \bm{\bar{\delta}}=(0,0,-1),\\
  &\ \ \,\text{(b)}:
  \ \tilde{f}_n(\bm{\lambda})=1,
  \ \ \tilde{b}_n(\bm{\lambda})=-(n+a-1)(n+b),
  \ \ \bm{\bar{\delta}}=(1,-1,0),\\
  &\ \ \,\text{(c)}:
  \ \tilde{f}_n(\bm{\lambda})=-(n+a)(n+b-1),
  \ \ \tilde{b}_n(\bm{\lambda})=1,
  \ \ \bm{\bar{\delta}}=(-1,1,0),\\
  &\ \ \,\text{(d)}:
  \ \tilde{f}_n(\bm{\lambda})
  =\mathcal{E}_N(\bm{\lambda})-\mathcal{E}_n(\bm{\lambda}),
  \ \ \tilde{b}_n(\bm{\lambda})=1,
  \ \ \bm{\bar{\delta}}=(0,0,1),\!\\
%
%%%%%%%%%%%%%%%%%%%%%%%%%%%%%%%%
% Krawtchouk (K)               %
%%%%%%%%%%%%%%%%%%%%%%%%%%%%%%%%
  \text{K}&:\text{(a)}:
  \ \tilde{f}_n(\bm{\lambda})=1,
  \ \ \tilde{b}_n(\bm{\lambda})
  =\mathcal{E}_{N+1}(\bm{\lambda})-\mathcal{E}_n(\bm{\lambda}),
  \ \ \bm{\bar{\delta}}=(0,-1),\\
  &\ \ \,\text{(b)}:
  \ \tilde{f}_n(\bm{\lambda})
  =\mathcal{E}_N(\bm{\lambda})-\mathcal{E}_n(\bm{\lambda}),
  \ \ \tilde{b}_n(\bm{\lambda})=1,
  \ \ \bm{\bar{\delta}}=(0,1),\\
%
%%%%%%%%%%%%%%%%%%%%%%%%%%%%%%%%
% Racah (R)                    %
%%%%%%%%%%%%%%%%%%%%%%%%%%%%%%%%
  \text{R}&:\text{(a)}:
  \,\tilde{f}_n(\bm{\lambda})=1,
  \ \,\tilde{b}_n(\bm{\lambda})
  =\mathcal{E}_{N+1}(\bm{\lambda})-\mathcal{E}_n(\bm{\lambda}),
  \ \bm{\bar{\delta}}=(1,0,0,1),\\
  &\ \ \,\text{(b)}:
  \,\tilde{f}_n(\bm{\lambda})=1,
  \ \,\tilde{b}_n(\bm{\lambda})=-(n+b-1)(n+c-d-N),
  \ \bm{\bar{\delta}}=(0,1,0,1),\\
  &\ \ \,\text{(c)}:
  \,\tilde{f}_n(\bm{\lambda})=1,
  \ \,\tilde{b}_n(\bm{\lambda})=-(n+c-1)(n+b-d-N),
  \ \bm{\bar{\delta}}=(0,0,1,1),\\
  &\ \ \,\text{(d)}:
  \,\tilde{f}_n(\bm{\lambda})=-(n+c)(n+b-d-N-1),
  \ \tilde{b}_n(\bm{\lambda})=1,
  \ \bm{\bar{\delta}}=(0,0,-1,-1),\!\!\\
  &\ \ \,\text{(e)}:
  \,\tilde{f}_n(\bm{\lambda})=-(n+b)(n+c-d-N-1),
  \ \tilde{b}_n(\bm{\lambda})=1,
  \ \bm{\bar{\delta}}=(0,-1,0,-1),\!\!\\
  &\ \ \,\text{(f)}:
  \,\tilde{f}_n(\bm{\lambda})
  =\mathcal{E}_N(\bm{\lambda})-\mathcal{E}_n(\bm{\lambda}),
  \ \tilde{b}_n(\bm{\lambda})=1,
  \ \bm{\bar{\delta}}=(-1,0,0,-1),\!\!\!\\
%
%%%%%%%%%%%%%%%%%%%%%%%%%%%%%%%%
% dual Hahn (dH)               %
%%%%%%%%%%%%%%%%%%%%%%%%%%%%%%%%
  \text{dH}&:\text{(a)}:
  \ \tilde{f}_n(\bm{\lambda})=1,
  \ \ \tilde{b}_n(\bm{\lambda})
  =\mathcal{E}_{N+1}(\bm{\lambda})-\mathcal{E}_n(\bm{\lambda}),
  \ \ \bm{\bar{\delta}}=(0,1,-1),\\
  &\ \ \,\text{(b)}:
  \ \tilde{f}_n(\bm{\lambda})=1,
  \ \ \tilde{b}_n(\bm{\lambda})=-(n+a-1),
  \ \ \bm{\bar{\delta}}=(1,0,0),\\
  &\ \ \,\text{(c)}:
  \ \tilde{f}_n(\bm{\lambda})=1,
  \ \ \tilde{b}_n(\bm{\lambda})=-(n-b-N+1),
  \ \ \bm{\bar{\delta}}=(0,1,0),\\
  &\ \ \,\text{(d)}:
  \ \tilde{f}_n(\bm{\lambda})=-(n-b-N),
  \ \ \tilde{b}_n(\bm{\lambda})=1,
  \ \ \bm{\bar{\delta}}=(0,-1,0),\\
  &\ \ \,\text{(e)}:
  \ \tilde{f}_n(\bm{\lambda})=-(n+a),
  \ \ \tilde{b}_n(\bm{\lambda})=1,
  \ \ \bm{\bar{\delta}}=(-1,0,0),\\
  &\ \ \,\text{(f)}:
  \ \tilde{f}_n(\bm{\lambda})
  =\mathcal{E}_N(\bm{\lambda})-\mathcal{E}_n(\bm{\lambda}),
  \ \ \tilde{b}_n(\bm{\lambda})=1,
  \ \ \bm{\bar{\delta}}=(0,-1,1),\\
%
%%%%%%%%%%%%%%%%%%%%%%%%%%%%%%%%%%%%%%
% dual quantum q-Krawtchouk (dq$q$K) %
%%%%%%%%%%%%%%%%%%%%%%%%%%%%%%%%%%%%%%
  \text{dq$q$K}&:\text{(a)}:
  \ \tilde{f}_n(\bm{\lambda})=1,
  \ \ \tilde{b}_n(\bm{\lambda})
  =\mathcal{E}_{N+1}(\bm{\lambda})-\mathcal{E}_n(\bm{\lambda}),
  \ \ \bm{\bar{\delta}}=(1,-1),\\
  &\ \ \,\text{(b)}:
  \ \tilde{f}_n(\bm{\lambda})=1,
  \ \ \tilde{b}_n(\bm{\lambda})=-q^{-n}(1-p^{-1}q^{n-N}),
  \ \ \bm{\bar{\delta}}=(1,0),\\
  &\ \ \,\text{(c)}:
  \ \tilde{f}_n(\bm{\lambda})=-q^{-n}(1-p^{-1}q^{n-N-1}),
  \ \ \tilde{b}_n(\bm{\lambda})=1,
  \ \ \bm{\bar{\delta}}=(-1,0),\\
  &\ \ \,\text{(d)}:
  \ \tilde{f}_n(\bm{\lambda})
  =\mathcal{E}_N(\bm{\lambda})-\mathcal{E}_n(\bm{\lambda}),
  \ \ \tilde{b}_n(\bm{\lambda})=1,
  \ \ \bm{\bar{\delta}}=(-1,1),\\
%
%%%%%%%%%%%%%%%%%%%%%%%%%%%%%%%%
% q-Hahn (qH)                  %
%%%%%%%%%%%%%%%%%%%%%%%%%%%%%%%%
  \text{$q$H}&:\text{(a)}:
  \ \tilde{f}_n(\bm{\lambda})=1,
  \ \ \tilde{b}_n(\bm{\lambda})
  =\mathcal{E}_{N+1}(\bm{\lambda})-\mathcal{E}_n(\bm{\lambda}),
  \ \ \bm{\bar{\delta}}=(0,0,-1),\\
  &\ \ \,\text{(b)}:
  \ \tilde{f}_n(\bm{\lambda})=1,
  \ \ \tilde{b}_n(\bm{\lambda})=-q^{-n}(1-aq^{n-1})(1-bq^n),
  \ \ \bm{\bar{\delta}}=(1,-1,0),\\
  &\ \ \,\text{(c)}:
  \ \tilde{f}_n(\bm{\lambda})=-q^{-n}(1-aq^n)(1-bq^{n-1}),
  \ \ \tilde{b}_n(\bm{\lambda})=1,
  \ \ \bm{\bar{\delta}}=(-1,1,0),\\
  &\ \ \,\text{(d)}:
  \ \tilde{f}_n(\bm{\lambda})
  =\mathcal{E}_N(\bm{\lambda})-\mathcal{E}_n(\bm{\lambda}),
  \ \ \tilde{b}_n(\bm{\lambda})=1,
  \ \ \bm{\bar{\delta}}=(0,0,1),\\
%
%%%%%%%%%%%%%%%%%%%%%%%%%%%%%%%%
% q-Krawtchouk ($q$K)          %
%%%%%%%%%%%%%%%%%%%%%%%%%%%%%%%%
  \text{$q$K}&:\text{(a)}:
  \ \tilde{f}_n(\bm{\lambda})=1,
  \ \ \tilde{b}_n(\bm{\lambda})
  =\mathcal{E}_{N+1}(\bm{\lambda})-\mathcal{E}_n(\bm{\lambda}),
  \ \ \bm{\bar{\delta}}=(0,-1),\\
  &\ \ \,\text{(b)}:
  \ \tilde{f}_n(\bm{\lambda})
  =\mathcal{E}_N(\bm{\lambda})-\mathcal{E}_n(\bm{\lambda}),
  \ \ \tilde{b}_n(\bm{\lambda})=1,
  \ \ \bm{\bar{\delta}}=(0,1),\\
%
%%%%%%%%%%%%%%%%%%%%%%%%%%%%%%%%
% quantum q-Krawtchouk (q$q$K) %
%%%%%%%%%%%%%%%%%%%%%%%%%%%%%%%%
  \text{q$q$K}&:\text{(a)}:
  \ \tilde{f}_n(\bm{\lambda})=1,
  \ \ \tilde{b}_n(\bm{\lambda})
  =\mathcal{E}_{N+1}(\bm{\lambda})-\mathcal{E}_n(\bm{\lambda}),
  \ \ \bm{\bar{\delta}}=(0,-1),\\
  &\ \ \,\text{(b)}:
  \ \tilde{f}_n(\bm{\lambda})=1,
  \ \ \tilde{b}_n(\bm{\lambda})=-p^{-1}q^{-1}(1-pq^{n+1}),
  \ \ \bm{\bar{\delta}}=(-1,0),\\
  &\ \ \,\text{(c)}:
  \ \tilde{f}_n(\bm{\lambda})=-p^{-1}(1-pq^n)
  \ \ \tilde{b}_n(\bm{\lambda})=1,
  \ \ \bm{\bar{\delta}}=(1,0),\\
  &\ \ \,\text{(d)}:
  \ \tilde{f}_n(\bm{\lambda})
  =\mathcal{E}_N(\bm{\lambda})-\mathcal{E}_n(\bm{\lambda}),
  \ \ \tilde{b}_n(\bm{\lambda})=1,
  \ \ \bm{\bar{\delta}}=(0,1),\\
%
%%%%%%%%%%%%%%%%%%%%%%%%%%%%%%%%
% affine q-Krawtchouk (a$q$K) %
%%%%%%%%%%%%%%%%%%%%%%%%%%%%%%%%
  \text{a$q$K}&:\text{(a)}:
  \ \tilde{f}_n(\bm{\lambda})=1,
  \ \ \tilde{b}_n(\bm{\lambda})
  =\mathcal{E}_{N+1}(\bm{\lambda})-\mathcal{E}_n(\bm{\lambda}),
  \ \ \bm{\bar{\delta}}=(0,-1),\\
  &\ \ \,\text{(b)}:
  \ \tilde{f}_n(\bm{\lambda})=1,
  \ \ \tilde{b}_n(\bm{\lambda})=-q^{-n}(1-pq^n),
  \ \ \bm{\bar{\delta}}=(1,0),\\
  &\ \ \,\text{(c)}:
  \ \tilde{f}_n(\bm{\lambda})=-q^{-n}(1-pq^{n+1}),
  \ \ \tilde{b}_n(\bm{\lambda})=1,
  \ \ \bm{\bar{\delta}}=(-1,0),\\
  &\ \ \,\text{(d)}:
  \ \tilde{f}_n(\bm{\lambda})
  =\mathcal{E}_N(\bm{\lambda})-\mathcal{E}_n(\bm{\lambda}),
  \ \ \tilde{b}_n(\bm{\lambda})=1,
  \ \ \bm{\bar{\delta}}=(0,1),\\
%
%%%%%%%%%%%%%%%%%%%%%%%%%%%%%%%%
% dual q-Hahn (dqH)            %
%%%%%%%%%%%%%%%%%%%%%%%%%%%%%%%%
  \text{d$q$H}&:\text{(a)}:
  \ \tilde{f}_n(\bm{\lambda})=1,
  \ \ \tilde{b}_n(\bm{\lambda})
  =\mathcal{E}_{N+1}(\bm{\lambda})-\mathcal{E}_n(\bm{\lambda}),
  \ \ \bm{\bar{\delta}}=(0,1,-1),\\
  &\ \ \,\text{(b)}:
  \ \tilde{f}_n(\bm{\lambda})=1,
  \ \ \tilde{b}_n(\bm{\lambda})=-q^{-n}(1-aq^{n-1}),
  \ \ \bm{\bar{\delta}}=(1,0,0),\\
  &\ \ \,\text{(c)}:
  \ \tilde{f}_n(\bm{\lambda})=1,
  \ \ \tilde{b}_n(\bm{\lambda})=-q^{-n}(1-b^{-1}q^{n-N+1}),
  \ \ \bm{\bar{\delta}}=(0,1,0),\\
  &\ \ \,\text{(d)}:
  \ \tilde{f}_n(\bm{\lambda})=-q^{-n}(1-b^{-1}q^{n-N}),
  \ \ \tilde{b}_n(\bm{\lambda})=1,
  \ \ \bm{\bar{\delta}}=(0,-1,0),\\
  &\ \ \,\text{(e)}:
  \ \tilde{f}_n(\bm{\lambda})=-q^{-n}(1-aq^n),
  \ \ \tilde{b}_n(\bm{\lambda})=1,
  \ \ \bm{\bar{\delta}}=(-1,0,0),\\
  &\ \ \,\text{(f)}:
  \ \tilde{f}_n(\bm{\lambda})
  =\mathcal{E}_N(\bm{\lambda})-\mathcal{E}_n(\bm{\lambda}),
  \ \ \tilde{b}_n(\bm{\lambda})=1,
  \ \ \bm{\bar{\delta}}=(0,-1,1),\\
%
%%%%%%%%%%%%%%%%%%%%%%%%%%%%%%%%
% dual q-Krawtchouk (dqK)      %
%%%%%%%%%%%%%%%%%%%%%%%%%%%%%%%%
  \text{d$q$K}&:\text{(a)}:
  \ \tilde{f}_n(\bm{\lambda})=1,
  \ \ \tilde{b}_n(\bm{\lambda})
  =\mathcal{E}_{N+1}(\bm{\lambda})-\mathcal{E}_n(\bm{\lambda}),
  \ \ \bm{\bar{\delta}}=(1,-1),\\
  &\ \ \,\text{(b)}:
  \ \tilde{f}_n(\bm{\lambda})=1,
  \ \ \tilde{b}_n(\bm{\lambda})=-q^{-n},
  \ \ \bm{\bar{\delta}}=(1,0),\\
  &\ \ \,\text{(c)}:
  \ \tilde{f}_n(\bm{\lambda})=-q^{-n},
  \ \ \tilde{b}_n(\bm{\lambda})=1,
  \ \ \bm{\bar{\delta}}=(-1,0),\\
  &\ \ \,\text{(d)}:
  \ \tilde{f}_n(\bm{\lambda})
  =\mathcal{E}_N(\bm{\lambda})-\mathcal{E}_n(\bm{\lambda}),
  \ \ \tilde{b}_n(\bm{\lambda})=1,
  \ \ \bm{\bar{\delta}}=(-1,1),\\
%
%%%%%%%%%%%%%%%%%%%%%%%%%%%%%%%%
% Meixner (M)                  %
%%%%%%%%%%%%%%%%%%%%%%%%%%%%%%%%
  \text{M}&:\text{(a)}:
  \ \tilde{f}_n(\bm{\lambda})=1,
  \ \ \tilde{b}_n(\bm{\lambda})=-(n+\beta-1),
  \ \ \bm{\bar{\delta}}=(1,0),\\
  &\ \ \,\text{(b)}:
  \ \tilde{f}_n(\bm{\lambda})=-(n+\beta),
  \ \ \tilde{b}_n(\bm{\lambda})=1,
  \ \ \bm{\bar{\delta}}=(-1,0),\\
%
%%%%%%%%%%%%%%%%%%%%%%%%%%%%%%%%
% little $q$-Jacobi (l$q$J)    %
%%%%%%%%%%%%%%%%%%%%%%%%%%%%%%%%
  \text{l$q$J}&:\text{(a)}:
  \ \tilde{f}_n(\bm{\lambda})=1,
  \ \ \tilde{b}_n(\bm{\lambda})=-q^{-n}(1-aq^n)(1-bq^{n-1}),
  \ \ \bm{\bar{\delta}}=(-1,1),\\
  &\ \ \,\text{(b)}:
  \ \tilde{f}_n(\bm{\lambda})=-q^{-n}(1-aq^{n-1})(1-bq^n),
  \ \ \tilde{b}_n(\bm{\lambda})=1,
  \ \ \bm{\bar{\delta}}=(1,-1),\\
%
%%%%%%%%%%%%%%%%%%%%%%%%%%%%%%%%
% little $q$-Laguerre (l$q$L)  %
%%%%%%%%%%%%%%%%%%%%%%%%%%%%%%%%
  \text{l$q$L}&:\text{(a)}:
  \ \tilde{f}_n(\bm{\lambda})=1,
  \ \ \tilde{b}_n(\bm{\lambda})=-q^{-n}(1-aq^n),
  \ \ \bm{\bar{\delta}}=-1,\\
  &\ \ \,\text{(b)}:
  \ \tilde{f}_n(\bm{\lambda})=-q^{-n}(1-aq^{n-1}),
  \ \ \tilde{b}_n(\bm{\lambda})=1,
  \ \ \bm{\bar{\delta}}=1,\\
%
%%%%%%%%%%%%%%%%%%%%%%%%%%%%%%%%
% $q$-Meixner ($q$M)           %
%%%%%%%%%%%%%%%%%%%%%%%%%%%%%%%%
  \text{$q$M}&:\text{(a)}:
  \ \tilde{f}_n(\bm{\lambda})=1,
  \ \ \tilde{b}_n(\bm{\lambda})=q^n+cq^{-1},
  \ \ \bm{\bar{\delta}}=(0,1),\\
  &\ \ \,\text{(b)}:
  \ \tilde{f}_n(\bm{\lambda})=1,
  \ \ \tilde{b}_n(\bm{\lambda})=-b^{-1}(1-bq^n),
  \ \ \bm{\bar{\delta}}=(1,0),\\
  &\ \ \,\text{(c)}:
  \ \tilde{f}_n(\bm{\lambda})=-b^{-1}q^{-1}(1-bq^{n+1}),
  \ \ \tilde{b}_n(\bm{\lambda})=1,
  \ \ \bm{\bar{\delta}}=(-1,0),\\
  &\ \ \,\text{(d)}:
  \ \tilde{f}_n(\bm{\lambda})=q^n+c,
  \ \ \tilde{b}_n(\bm{\lambda})=1,
  \ \ \bm{\bar{\delta}}=(0,-1),\\
%
%%%%%%%%%%%%%%%%%%%%%%%%%%%%%%%%%%%%%
% Al-Salam-Carlitz $\II$ (ASC$\II$) %
%%%%%%%%%%%%%%%%%%%%%%%%%%%%%%%%%%%%%
  \text{ASC$\II$}&:\text{(a)}:
  \ \tilde{f}_n(\bm{\lambda})=1,
  \ \ \tilde{b}_n(\bm{\lambda})=q^n,
  \ \ \bm{\bar{\delta}}=1,\\
  &\ \ \,\text{(b)}:
  \ \tilde{f}_n(\bm{\lambda})=q^n,
  \ \ \tilde{b}_n(\bm{\lambda})=1,
  \ \ \bm{\bar{\delta}}=-1,\\
%
%%%%%%%%%%%%%%%%%%%%%%%%%%%%%%%%
% $q$-Charlier ($q$C)          %
%%%%%%%%%%%%%%%%%%%%%%%%%%%%%%%%
  \text{$q$C}&:\text{(a)}:
  \ \tilde{f}_n(\bm{\lambda})=1,
  \ \ \tilde{b}_n(\bm{\lambda})=q^n+aq^{-1},
  \ \ \bm{\bar{\delta}}=1,\\
  &\ \ \,\text{(b)}:
  \ \tilde{f}_n(\bm{\lambda})=q^n+a,
  \ \ \tilde{b}_n(\bm{\lambda})=1,
  \ \ \bm{\bar{\delta}}=-1.
\end{align}

%%%%%%%%%%%%%%%%%%%%%%%%%%%%%%%%%%%%%%%%%%%%%%%%%%%%%
%                                                   %
% A.4 Data for \S\,\ref{sec:newrdQMJ}               %
%                                                   %
%%%%%%%%%%%%%%%%%%%%%%%%%%%%%%%%%%%%%%%%%%%%%%%%%%%%%
\subsection{Data for \S\,\ref{sec:newrdQMJ}}
\label{app:rdQMJ}

We present explicit forms of $B^{\text{J}}_1(\eta)$, $D^{\text{J}}_1(\eta)$,
$\tilde{f}_n$, $\tilde{b}_n$ and $\bm{\bar{\delta}}$ in \S\,\ref{sec:newrdQMJ}.
The potential functions $B^{\text{J}}_2(\eta)$ and $D^{\text{J}}_2(\eta)$ can
be obtained from \eqref{BJ=BJ1BJ2}.

\noindent
Potential functions $B^{\text{J}}_1(\eta;\bm{\lambda})$ and
$D^{\text{J}}_1(\eta;\bm{\lambda})$:
\begin{align}
%%%%%%%%%%%%%%%%%%%%%%%%%%%%%%%%
% big $q$-Laguerre (b$q$L)     %
%%%%%%%%%%%%%%%%%%%%%%%%%%%%%%%%
  \text{b$q$L}&:\text{(a)}:
  \ B^{\text{J}}_1(\eta;\bm{\lambda})=\frac{\eta^{-1}a(1-\eta)}{1-a},
  \ \ D^{\text{J}}_1(\eta;\bm{\lambda})=\frac{\eta^{-1}(aq-\eta)}{a-1},\\
  &\ \ \,\text{(b)}:
  \ B^{\text{J}}_1(\eta;\bm{\lambda})=\frac{\eta^{-1}b(1-\eta)}{1-b},
  \ \ D^{\text{J}}_1(\eta;\bm{\lambda})=\frac{\eta^{-1}(\eta-bq)}{1-b},\\
  &\ \ \,\text{(c)}:
  \ B^{\text{J}}_1(\eta;\bm{\lambda})=(bq-1)\eta^{-1}a,
  \ \ D^{\text{J}}_1(\eta;\bm{\lambda})=(1-bq)\eta^{-1}(aq-\eta),\\
  &\ \ \,\text{(d)}:
  \ B^{\text{J}}_1(\eta;\bm{\lambda})=(aq-1)\eta^{-1}b,
  \ \ D^{\text{J}}_1(\eta;\bm{\lambda})=(aq-1)\eta^{-1}(\eta-bq),\\
%
%%%%%%%%%%%%%%%%%%%%%%%%%%%%%%%%%%%
% Al-Salam-Carlitz $\I$ (ASC$\I$) %
%%%%%%%%%%%%%%%%%%%%%%%%%%%%%%%%%%%
  \text{ASC$\I$}&:\text{(a)}:
  \ B^{\text{J}}_1(\eta;\bm{\lambda})=\eta^{-1}q^{-1},
  \ \ D^{\text{J}}_1(\eta;\bm{\lambda})=-\eta^{-1}(1-\eta),\\
  &\ \ \,\text{(b)}:
  \ B^{\text{J}}_1(\eta;\bm{\lambda})=-\eta^{-1}aq^{-1},
  \ \ D^{\text{J}}_1(\eta;\bm{\lambda})=-\eta^{-1}(\eta-a),\\
%
%%%%%%%%%%%%%%%%%%%%%%%%%%%%%%%%
% $q$-Laguerre ($q$L)          %
%%%%%%%%%%%%%%%%%%%%%%%%%%%%%%%%
  \text{$q$L}&:\text{(a)}:
  \ B^{\text{J}}_1(\eta;\bm{\lambda})=\eta^{-1}(1+\eta),
  \ \ D^{\text{J}}_1(\eta;\bm{\lambda})=-\eta^{-1}q,\\
  &\ \ \,\text{(b)}:
  \ B^{\text{J}}_1(\eta;\bm{\lambda})=1,
  \ \ D^{\text{J}}_1(\eta;\bm{\lambda})=-a^{-1}.
\end{align}
Constants $\tilde{f}^{\text{J}}_n(\bm{\lambda})$,
$\tilde{b}^{\text{J}}_n(\bm{\lambda})$ and $\bm{\bar{\delta}}$:
\begin{align}
%%%%%%%%%%%%%%%%%%%%%%%%%%%%%%%%
% big $q$-Laguerre (b$q$L)     %
%%%%%%%%%%%%%%%%%%%%%%%%%%%%%%%%
  \text{b$q$L}&:\text{(a)}:
  \ \tilde{f}^{\text{J}}_n(\bm{\lambda})=1,
  \ \ \tilde{b}^{\text{J}}_n(\bm{\lambda})=-q^{-n}(1-aq^n),
  \ \ \bm{\bar{\delta}}=(1,0),\\
  &\ \ \,\text{(b)}:
  \ \tilde{f}^{\text{J}}_n(\bm{\lambda})=1,
  \ \ \tilde{b}^{\text{J}}_n(\bm{\lambda})=-q^{-n}(1-bq^n),
  \ \ \bm{\bar{\delta}}=(0,1),\\
  &\ \ \,\text{(c)}:
  \ \tilde{f}^{\text{J}}_n(\bm{\lambda})=-q^{-n}(1-bq^{n+1}),
  \ \ \tilde{b}^{\text{J}}_n(\bm{\lambda})=1,
  \ \ \bm{\bar{\delta}}=(0,-1),\\
  &\ \ \,\text{(d)}:
  \ \tilde{f}^{\text{J}}_n(\bm{\lambda})=-q^{-n}(1-aq^{n+1}),
  \ \ \tilde{b}^{\text{J}}_n(\bm{\lambda})=1,
  \ \ \bm{\bar{\delta}}=(-1,0),\\
%
%%%%%%%%%%%%%%%%%%%%%%%%%%%%%%%%%%%
% Al-Salam-Carlitz $\I$ (ASC$\I$) %
%%%%%%%%%%%%%%%%%%%%%%%%%%%%%%%%%%%
  \text{ASC$\I$}&:\text{(a)}:
  \ \tilde{f}^{\text{J}}_n(\bm{\lambda})=1,
  \ \ \tilde{b}^{\text{J}}_n(\bm{\lambda})=-q^{-n},
  \ \ \bm{\bar{\delta}}=-1,\\
  &\ \ \,\text{(b)}:
  \ \tilde{f}^{\text{J}}_n(\bm{\lambda})=-q^{-n},
  \ \ \tilde{b}^{\text{J}}_n(\bm{\lambda})=1,
  \ \ \bm{\bar{\delta}}=1,\\
%
%%%%%%%%%%%%%%%%%%%%%%%%%%%%%%%%
% $q$-Laguerre ($q$L)          %
%%%%%%%%%%%%%%%%%%%%%%%%%%%%%%%%
  \text{$q$L}&:\text{(a)}:
  \ \tilde{f}^{\text{J}}_n(\bm{\lambda})=1,
  \ \ \tilde{b}^{\text{J}}_n(\bm{\lambda})=-a^{-1}q^{-1}(1-aq^{n+1}),
  \ \ \bm{\bar{\delta}}=-1,\\
  &\ \ \,\text{(b)}:
  \ \tilde{f}^{\text{J}}_n(\bm{\lambda})=-a^{-1}(1-aq^n),
  \ \ \tilde{b}^{\text{J}}_n(\bm{\lambda})=1,
  \ \ \bm{\bar{\delta}}=1.
\end{align}

%%%%%%%%%%%%%%%%%%%%%%%%%%%%%%%%%%%%%%%%%%%%%%%%%%%%%%%%%%%%%%%
%                                                             %
%  References                                                 %
%                                                             %
%%%%%%%%%%%%%%%%%%%%%%%%%%%%%%%%%%%%%%%%%%%%%%%%%%%%%%%%%%%%%%%

\end{document}